\newif\ifarxiv
\arxivtrue 

\ifarxiv 
  \documentclass[11pt,a4paper]{article}
  \usepackage[margin=2cm]{geometry}
  \usepackage{cuted}
  \usepackage{amsmath,amssymb,amsthm}
  \usepackage[affil-it]{authblk}
  \usepackage[colorlinks,bookmarksopen,bookmarksnumbered,citecolor=red,urlcolor=red]{hyperref}  
  \newcommand{\sep}{, }
  \newcommand{\PACS}{\noindent\textbf{PACS:} }
  \newcommand{\MSC}{\noindent\textbf{2000 MSC:} }
  \newcommand{\appref}[1]{Appendix~\ref{#1}}
\else 
   \documentclass[preprint,11pt]{elsarticle}
   \usepackage[margin=2cm]{geometry}
   \newcommand{\appref}[1]{\ref{#1}}
\fi 

\usepackage{amssymb}
\usepackage{amsmath}
\usepackage{xfrac}
\usepackage{soul}
\usepackage{color}
\usepackage{multirow}
\usepackage{mathrsfs}
\usepackage{hyperref}
\hypersetup{
    colorlinks=true,
    linkcolor=blue,
    citecolor=red,
    filecolor=magenta,      
    urlcolor=blue
    }
\usepackage{algorithm}
\usepackage{algpseudocode}
\usepackage{array}
\newcolumntype{x}[1]{>{\centering\arraybackslash\hspace{0pt}}p{#1}}

\newcommand{\jump}[1]{[\![ #1]\!]}
\newcommand{\avg}[1]{\{\!\{#1\}\!\}}
\newcommand{\Id}{\operatorname{Id}}
\newcommand{\impl}{{(\text{im})}}
\newcommand{\expl}{{(\text{ex})}}
\newcommand{\figdir}{./}

\ifarxiv 
  \author[*]{Eike~Hermann~M\"{u}ller}
  \affil[ ]{Department of Mathematical Sciences, University of Bath, Bath BA2 7AY, Bath, United Kingdom}
  \affil[*]{Email: \texttt{e.mueller@bath.ac.uk}}
\else 
  \journal{Journal of Computational Physics}
\fi 
\title{A Higher-order Hybridisable Discontinuous Galerkin IMEX method for the incompressible Euler equations}
\begin{document}

\ifarxiv 
    \maketitle
\else 
    \begin{frontmatter}


        \author{Eike Hermann Mueller\corref{cor1}}
        \ead{e.mueller@bath.ac.uk}
        \cortext[cor1]{Corresponding author}
        \affiliation{organization={Department of Mathematical Sciences, University of Bath},
            addressline={Claverton Down},
            city={Bath},
            postcode={BA2 7AY},
            country={United Kingdom}}
        \fi 

        \begin{abstract}
            The incompressible Euler equations are an important model system in computational fluid dynamics. Fast high-order methods for the solution of this time-dependent system of partial differential equations are of particular interest: due to their exponential convergence in the polynomial degree they can make efficient use of computational resources. To address this challenge we describe a novel timestepping method which combines a hybridised Discontinuous Galerkin method for the spatial discretisation with IMEX timestepping schemes, thus achieving high-order accuracy in both space and time. The computational bottleneck is the solution of a (block-) sparse linear system to compute updates to pressure and velocity at each stage of the IMEX integrator. Following Chorin's projection approach, this update of the velocity and pressure fields is split into two stages. As a result, the hybridised equation for the implicit pressure-velocity problem is reduced to the well-known system which arises in hybridised mixed formulations of the Poisson- or diffusion problem and for which efficient multigrid preconditioners have been developed. Splitting errors can be reduced systematically by embedding this update into a preconditioned Richardson iteration. The accuracy and efficiency of the new method is demonstrated numerically for two time-dependent testcases that have been previously studied in the literature.
        \end{abstract}


        \ifarxiv 
        \noindent\textbf{Keywords:}
        \else 
        \begin{keyword}
            \fi 
            Computational Fluid Dynamics\sep Incompressible Euler equations\sep IMEX timestepping\sep Hybridisable Discontinuous Galerkin\sep Multigrid
            \ifarxiv 
            \else 
            \PACS 02.10.Ud\sep 02.30.Jr\sep 02.60.Lj\sep 02.60.-x\sep 02.70.Dh\sep 47.10.ab\sep 47.11.-j\sep 47.11.Fg
            \MSC 35J55\sep 35L65\sep 35Q35\sep 65L06\sep 65M55\sep 65M60\sep 65N30\sep 65N55
        \end{keyword}
    \end{frontmatter}
\fi 


\section{Introduction}\label{sec:introduction}
We consider the incompressible Euler equations for velocity $Q(x,t)$ and pressure $p(x,t)$
\begin{subequations}
    \begin{align}
        \partial_t Q + (Q\cdot \nabla) Q + \nabla p & = f, \qquad\text{(momentum conservation)}\label{eqn:incompressible_euler_momentum}         \\
        \nabla\cdot Q                               & = 0 \qquad\text{(incompressibility constraint)}\label{eqn:incompressible_euler_constraint}
    \end{align}
\end{subequations}
in the $d$-dimensional domain $\Omega\subset \mathbb{R}^d$ with the initial condition $Q(x,0) = Q_0(x)$ for all $x\in\Omega$ and the boundary condition \mbox{$n\cdot Q\vert_{\partial \Omega}=0$} on $\partial\Omega$. In the momentum conservation equation \eqref{eqn:incompressible_euler_momentum} $f=f(t)$ is a time dependent forcing function. It is worth reminding the reader that -- in contrast to its role in the \textit{compressible} equations of fluid flow -- the pressure acts as a Lagrange multiplier which enforces the incompressibility constraint in \eqref{eqn:incompressible_euler_constraint}. Unlike the velocity $Q$, which is a dynamic (or prognostic) variable, the pressure is purely diagnostic and can be recovered from $Q$ at any given time by solving the Poisson problem $-\Delta p = \nabla\cdot((Q\cdot\nabla) Q)- \nabla \cdot f$ that is obtained by taking the divergence of \eqref{eqn:incompressible_euler_momentum} while exploiting the incompressibility constraint in \eqref{eqn:incompressible_euler_constraint}.
\subsection{Numerical methods for the incompressible Euler- and Navier-Stokes equations}
In \cite{Guzman2016} two different spatial discretisations for the incompressible Euler equations in \eqref{eqn:incompressible_euler_momentum}, \eqref{eqn:incompressible_euler_constraint} are considered. The first method uses a H(div)-conforming velocity field, which is discretised in a Raviart Thomas space \cite{RaviartThomas1977}, and a Discontinuous Galerkin (DG) \cite{Reed1973} space for the pressure. The second method uses a DG discretisation for both pressure and velocity and requires stabilisation of the advection term. While both methods are of arbitrary order in space, only a simple fully implicit timestepper is used in \cite{Guzman2016}. This leads to two limitations: firstly, the temporal discretisation is of lowest order and hence ultimately a very small timestep size will be required to ensure that the temporal discretisation error balances the spatial error. Secondly, the linear system which has to be solved in each timestep is non-trivial. In fact, no further effort is made in \cite{Guzman2016} to make the solution of this system efficient since this is not the focus of the paper. In \cite{Ueckermann2016} these limitations are addressed, but the incompressible Navier Stokes equations are considered instead. The authors of \cite{Ueckermann2016} propose a hybridised Discontinuous Galerkin (HDG) method \cite{Cockburn2009,Cockburn2010} for the spatial discretisation, employ an implicit-explicit IMEX timestepping scheme \cite{Ascher1997,pareschi2005implicit} of arbitrary order and use Chorin's projection method \cite{Chorin1967,Chorin1968} to simplify the solution of the resulting linear systems.

The purpose of the present paper is to apply the ideas in \cite{Ueckermann2016} to the incompressible Euler equations in  \eqref{eqn:incompressible_euler_momentum}, \eqref{eqn:incompressible_euler_constraint} to obtain a computationally efficient high-order discretisation in space and time also for this fluid flow problem. This is non-trivial since the two sets of equations have fundamental differences. In particular -- in contrast to the Navier Stokes equations -- the momentum conservation equation in \eqref{eqn:incompressible_euler_momentum} does not contain any diffusion terms which usually help to stabilise numerical schemes.
\subsection{Efficient multigrid preconditioners for HDG discretisations}
The computational bottleneck of IMEX schemes is the solution of a system of linear equation at each stage. It is therefore crucial to use optimal solvers and preconditioners for this task. As will be discussed below, in our approach the problem reduces to the solution of a coupled system of linear equations which has the same structure as the mixed formulation of the diffusion problem considered in \cite{Cockburn2014}: for a given positive definite $d\times d$ matrix function $a(x)$ find a vector-valued field $q(x)$ and a scalar-valued field $u(x)$ such that $q+a(x)\nabla u = 0$, $\nabla \cdot q=0$ in the domain $\Omega$ with Dirichlet boundary conditions on $\partial \Omega$. The authors of \cite{Cockburn2014} propose to solve this problem with a HDG method \cite{Cockburn2009,Cockburn2010}: for this, the numerical fluxes are expressed in terms of a trace-field $\lambda$ that is defined on a DG space on the set of mesh-facets. Elimination of $q$ and $u$ leads to a significantly smaller system for $\lambda$, which is then solved with a non-nested multigrid algorithm. A simple block-Jacobi method is employed to smooth the trace field on the finest level. The next-coarser level consists of piecewise linear, continuous functions on the same mesh, for this space a standard multigrid hierarchy can be constructed via $h$-coarsening. As shown in \cite{Cockburn2014}, the resulting algorithm shows optimal convergence and the cost of one multigrid cycle is proportional to the number of unknowns on the finest level. Here we use the same approach to construct efficient solvers for the sparse linear systems that arise in the IMEX method when it is applied to the incompressible Euler equations.
\subsection{Novelty and main results}
The key achievement reported in the present paper is the construction of an efficient method for the solution of the incompressible Euler equations based on the ideas outlined above. More specifically:
\begin{enumerate}
    \item We extend the spatial discretisation in \cite{Guzman2016} to a hybridisable DG method \cite{Cockburn2009,Cockburn2010}.
    \item We embed this spatial HDG discretisation into an IMEX \cite{Ascher1997,pareschi2005implicit} timestepping scheme.
    \item We use the projection method \cite{Chorin1967,Chorin1968} to reduce the linear systems that need to be solved to a form which can be tackled with known non-nested multigrid precondititioners \cite{Cockburn2014}.
    \item We demonstrate the high-order convergence of our method for a Taylor Green vortex \cite{Taylor1937} with known analytical solution and report results for a shear-flow problem \cite{Bell1989} also considered in \cite{Guzman2016}.
    \item We show that the method is efficient in the sense that empirically the cost per timestep grows in proportion to the problem size.
\end{enumerate}
All algorithms were implemented in the Firedrake \cite{FiredrakeUserManual} framework, which allows the straightforward high-level specification of the problem in weak form and supports composable solvers \cite{Kirby2017}. It should be stressed that in contrast to \cite{Guzman2016} no attempt is made here to analyse the numerical scheme theoretically. This is a significant endeavour which is left for a future publication.
\subsection{Related work}
DG formulations of the incompressible Navier Stokes equations can be found in standard textbooks such as \cite{hesthaven2007nodal,di2011mathematical} (however, the incompressible Euler equations are not discussed there). In \cite{hesthaven2007nodal} a splitting scheme is used for the time integration: in the second stage of this scheme, the velocity is projected to the space of divergence-free functions. The authors of \cite{liu2000high} develop a high-order discretisation for the two-dimensional incompressible Euler equations in the stream function formulation; the Navier Stokes equations in the high-Reynolds number region are also considered. Similar to the problem considered in our work, this approach requires the solution of a Poisson equation for the vorticity, which is discretised in a continuous finite element space, while a DG discretisation is used for the velocity. The momentum equation is treated fully explicitly with a low-order timestepping scheme. Hence, in contrast to the work presented here, only the spatial discretisation is high-order accurate. The focus of \cite{chen2021some} is the conservation of physical quantities in the solution of the incompressible Euler- and Navier-Stokes equations discretised with continuous and discontinuous Galerkin schemes. Again, in contrast to our work, a lower-order Crank-Nicolson timestepping scheme is used to advance the solution forward in time.
\subsubsection{Space-time discretisations}
Instead of combining a highly accurate spatial discretisation with higher-order timesteppers, some authors have proposed to discretise the problem simultaneously in space and time. By using suitable finite elements, this allows the systematic control of both the spatial and temporal discretisation error. The downside of this approach is that the (non-)linear systems that have to be solved are significantly more complicated and difficult to precondition. The authors of \cite{pesch2008discontinuous} describe a space-time DG discretisation of the incompressible (and compressible) Euler equations in $d$ dimensions. This approach allows high-order discretisations in space- and time, but requires the solution of a $d+1$ dimensional space-time problem. In \cite{pesch2008discontinuous} this is addressed by obtaining the steady-state solution of a pseudo-time integration problem. Related to the work in \cite{pesch2008discontinuous}, the authors of \cite{fambri2020discontinuous} consider space-time discretisations in the ADER-DG framework (see e.g. \cite{dumbser2008unified}) for the compressible Navier Stokes equations. ADER-DG methods combine a high-order explicit space-time predictor, which can be computed locally, with a corrector that uses the resulting fluxes to solve a global problem. While this results in a computationally very efficient scheme which is high-order in space and time and naturally supports adaptive mesh refinement, it is not obvious how the method can be extended to the \textit{incompressible} Euler equations, for which the incompressibility constraint has to be incorporated. In \cite{fambri2020discontinuous} a semi-implicit DG space-time discretisation of the incompressible Navier Stokes equations in constructed. For this, the nonlinear convection-diffusion terms are treated explicitly, while the incompressibility condition and the pressure gradients are incorporated implicitly. A space-time DG discretisation of the incompressible Navier-Stokes equations is also considered in \cite{rhebergen2013space}.
\subsubsection{Hybridisable DG methods}
Since they promise computational advantages, HDG discretisations \cite{Cockburn2009,Cockburn2010} have attracted significant interest in computational fluid dynamics. The authors of \cite{peraire2010hybridizable} develop a hybridised formulation for the compressible Navier-Stokes and Euler equations. Similar to the approach in the present paper, they also solve the resulting linear systems by eliminating the original state variables in favour of the trace-variables on the mesh skeleton. However, there is no detailled discussion of the efficient solution of this system or the development of higher-order timesteppers. To improve efficiency further, the authors of \cite{peraire2011embedded} use an embedded DG (EDG) method by aprroximating the numerical flux variables in a lower-order space. While this reduces the size of the hybridised system that has to be solved after elimination of the primal variables, this comes at the cost of reducing the order of convergence of some variables. Furthermore, as in \cite{peraire2010hybridizable}, only compressible flows are considered. While usually only the spatial discretisation is hybridised, the approach has also been extended to space-time formulations of the incompressible Navier-Stokes equations in \cite{rhebergen2012space}. While compared to native DG approaches in \cite{pesch2008discontinuous,fambri2020discontinuous,rhebergen2013space} the resulting system of equations is smaller, this approach still requires the solution of a complicated implicit space-time problem.
\subsubsection{DG methods for Compressible flows}
While not the topic of this paper, for completeness we also briefly review some related work on the compressible Euler- and Navier-Stokes equations. \cite{hartmann2006discontinuous} contains a detailed description of a DG discretisation of the compressible Euler- and Navier Stokes equations based on an interior penalty fomulation. However, only stationary solutions are considered. For the Euler problem they study the high-order convergence both theoretically and numerically for the Ringleb model \cite{chiocchia1984exact} for which the exact solution is known. Adaptive mesh refinement to capture shocks in the solution is also considered. The authors of \cite{williams2019analysis} analyse the stability of several DG discretisations of the compressible Euler equations; the methods they consider differ in the dynamical variables that are used; a fourth-order Singly Diagonally Implicit Runge–Kutta (SDIRK) \cite{alexander1977diagonally} is used to advance the solution in time. In \cite{zeifang2018efficient} a high-order DG scheme is combined with IMEX timestepping methods. For this, the flux is split into a non-stiff part which can be treated explicitly, and a stiff part, which requires the solution of an implicit problem. However, in contrast to the work presented here, the authors of \cite{zeifang2018efficient} consider compressible flows. DG discretisations of the compressible Navier Stokes equations are also discussed in \cite{jung2024behavior,qin2013discontinuous}.
\paragraph{Structure}
The rest of this paper is organised as follows: Section \ref{sec:methods} contains a detailed derivation of the numerical methods that are employed to discretise the incompressible Euler equations in space and advance them forward in time. The resulting algorithms are used to solve two well-known model problems in Section \ref{sec:results}, which also contains details on the implementation. We conclude in Section \ref{sec:conclusion} and some more technical details and additional results are relegated to the appendices.
\section{Methods}\label{sec:methods}
Our new timestepping method for the incompressible Euler equations is derived in four steps:
\begin{enumerate}
    \item We start by writing down the fully implicit DG method given in \cite{Guzman2016} in Section \ref{sec:dg_discretisation}.
    \item Next, we hybridise the method by introducing suitable fluxes in Section \ref{sec:hybridisation}.
    \item The splitting method \cite{Chorin1967,Chorin1968} which allows the efficient solution of the resulting linear systems is described in Section \ref{sec:splitting_method}.
    \item Finally, the full implicit-explicit (IMEX) integrator which combines all building blocks is written down in Section \ref{sec:imex}.
\end{enumerate}
\subsection{Notation}
\begin{figure}
    \begin{center}
        \includegraphics[width=0.3\linewidth]{\figdir/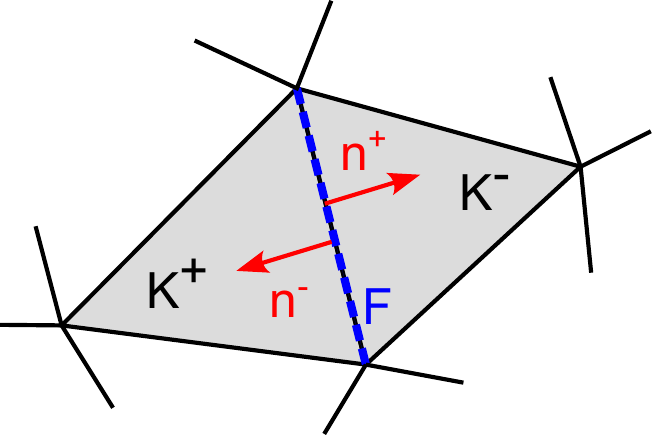}
        \caption{Two neighbouring cells $K^+$, $K^-$ with corresponding outward normals $n^+$, $n^-$.}
        \label{fig:two_cells}
    \end{center}
\end{figure}

Let $\Omega_h$ be a mesh covering the domain $\Omega$ where $h:=\max_{K\in\Omega}\{\operatorname{diam}(K)\}$ is the grid spacing. We write $\mathcal{E}_h$ for the skeleton of this mesh, i.e. the set of all facets. The set of interior facets is denoted by $\mathcal{E}^i_h$ and the set of boundary facets by $\mathcal{E}^\partial_h$ such that $\mathcal{E}_h= \mathcal{E}^i_h\cup \mathcal{E}^\partial_h$. For each interior facet $F\in\mathcal{E}_h^i$ we label the two cells touching this facet as $K^+$, $K^-$ with outward normals $n^+$, $n^-=-n^+$ respectively, see Fig.~\ref{fig:two_cells}.
Quantities associated with the two cells are identified with super- or subscripts $+$ and $-$. In particular we define
\begin{equation}
    a^\pm(x) := \lim_{\varepsilon\rightarrow 0^+} a(x-\varepsilon n^\pm)
    \qquad\text{for $x\in F$.}
\end{equation}
For a single cell $K$ we denote its outward normal by $n$ and write
\begin{equation}
    a^\pm(x) := \lim_{\varepsilon\rightarrow 0^+} a(x\mp \varepsilon n)
    \qquad\text{for $x\in \partial K$.}
\end{equation}
Further, for scalar quantities $a$ and vector-valued quantities $b$ we define the average $\avg{\cdot}$ and jump $\jump{\cdot}$ as
\begin{xalignat}{3}
    \avg{a} &:= \frac{1}{2}(a^++a^-), &
    \avg{b} &:= \frac{1}{2}(b^++b^-), &
    \jump{b\cdot n} &:= b^+\cdot n^+ + b^-\cdot n^-
\end{xalignat}
on interior facets $F\in\mathcal{E}_h^i$ and
\begin{xalignat}{3}
    \avg{a} &:= a, &
    \avg{b} &:= b, &
    \jump{b\cdot n} &:= b\cdot n
\end{xalignat}
on boundary facets $F\in\mathcal{E}_h^\partial$. The inner product of two tensors $A$, $B$ is given by
\begin{equation}
    A:B := \sum_{ij} A_{ij} B_{ji} = \operatorname{trace}(AB^\top).
\end{equation}
We use the shorthand
\begin{equation}
    (f)_A = \int_A f\;dx\qquad\text{for $A \subset \Omega$}
\end{equation}
for integrals over $d$-dimensional domains $A\subset \mathbb{R}^d$ (typically grid cells) and
\begin{equation}
    \langle f\rangle_S = \int_S f\;ds\qquad\text{for $S \subset \mathcal{E}_h$}
\end{equation}
for integrals over $d-1$-dimensional domains (typically facets or boundaries of grid cells). We write $h_F:=\langle 1\rangle_F$ for the area of a facet $F$.
\subsection{Discontinuous Galerkin spaces}
Let $P_k(A)$ be the space of polynomials of degree $k$ on $A$. With this we define the discontinuous Galerkin space as
\begin{equation}
    \text{DG}_{h,k} := \{ v\in L^2(\Omega) : v|_K \in P_k(K) \;\text{for all}\;K\in \Omega_h \}.
\end{equation}
The functions in this space are polynomials of degree $k$ in each grid cell, but no continuity is assumed across facets that separate neighbouring cells. In the rest of the paper we assume that the velocity $Q$ and pressure $p$ live in the following spaces, which are also used in \cite{Guzman2016}:
\begin{equation}
    \begin{aligned}
        Q \in V_Q & := [\text{DG}_{h,k+1}]^d,                    \\
        p \in V_p & := \text{DG}_{h,k}\cap \{v:(v)_{\Omega}=0\}. \\
    \end{aligned}\label{eqn:function_spaces}
\end{equation}
Recall that only the gradient of the pressure appears in \eqref{eqn:incompressible_euler_momentum}, \eqref{eqn:incompressible_euler_constraint} and the strong constraint $(p)_\Omega=\int_\Omega p\;dx=0$ on the pressure guarantees that the solution is unique.
\subsection{Original fully-implicit DG method}\label{sec:dg_discretisation}
For reference, we start by writing down the DG-based fully implicit timestepping method in \cite{Guzman2016}
(see (2.27) and (2.28) there), generalised to a non-zero right-hand side $f=f(t)$ for the momentum conservation equation in \eqref{eqn:incompressible_euler_momentum} and with a slightly different treatment of the no-flux boundary condition $n\cdot Q|_{\partial\Omega}=0$. For a given forcing function $f=f(t)\in V_Q$ we aim to find $(Q,p)\in V_Q\times V_p$ such that
\begin{equation}
    \begin{aligned}
        (\partial_t Q\cdot w)_{K} - (Q\otimes Q^\star:\nabla w)_{K} - (p\nabla\cdot w)_{K} + \langle \widehat{\sigma}: n\otimes w \rangle_{\partial K} & = (f\cdot w)_{K}, \\
        - (Q\cdot \nabla \psi)_{K} + \langle (\widehat{Q}\cdot n)\psi\rangle_{\partial K}                                                              & =0
    \end{aligned}\label{eqn:original_scheme}
\end{equation}
for all test functions $(w,\psi)\in V_Q\times V_p$ and all cells $K$. Here, $Q^\star =:\mathcal{P}(Q)$ is the interpolation of $Q$ onto the $\text{BDM}^0_{k+1}(\Omega_h)$ space as defined by (2.25) and (2.26) in \cite{Guzman2016}. The numerical fluxes $\widehat{\sigma}$ and $\widehat{Q}$ are given by
\begin{equation}
    \begin{aligned}
        \widehat{\sigma}   & := \avg{p} \Id  + \alpha h_F^{-1} \jump{Q\cdot n} \Id       + \begin{cases}
                                                                                               \avg{Q}\otimes Q^\star       & \text{for the central flux}, \\
                                                                                               Q^{\text{up}}\otimes Q^\star & \text{for the upwind flux}
                                                                                           \end{cases} \\
        \widehat{Q}\cdot n & := \avg{Q}\cdot n,
    \end{aligned}\label{eqn:flux_original}
\end{equation}
for some positive constant $\alpha>0$ and with the upwind velocity defined as
\begin{equation}
    Q^\text{up} = \begin{cases}
        Q^+ & \text{if $Q\cdot n\ge 0$} \\
        Q^- & \text{otherwise.}
    \end{cases}
    \label{eqn:upwind_velocity}
\end{equation}
Using a backward-Euler time discretisation but treating the forcing term on the right-hand side of \eqref{eqn:original_scheme} explicitly, this leads to the following time-stepping scheme for computing $(Q^{n+1},p^{n+1})\in V_Q\times V_p$ from $Q^n$ at timestep $n$ and time $t^n=n\Delta t$:
\begin{equation}
    \begin{aligned}
        (w\cdot Q^{n+1})_{\Omega_h} - \Delta t \left(f^\impl(w,Q^{n+1},\mathcal{P}(Q^n)) + \widetilde{g}(w,p^{n+1})\right) & = (w\cdot Q^{n})_{\Omega_h} + \Delta t f^\expl(w, Q^n,t^n) \\
        \text{subject to}\qquad\widetilde{\Gamma}(\psi,Q^{n+1})                                                            & = 0
    \end{aligned}\label{eqn:original_DG_scheme}
\end{equation}
for all test functions $(w,\psi) \in V_Q\times V_p$. Explicit expressions for the terms $f^\impl$, $f^\expl$ in \eqref{eqn:original_DG_scheme} are given by
\begin{subequations}
    \begin{align}
        f^\impl(w,Q,Q^\star) & := -\left(w\cdot (Q^\star\cdot \nabla) Q\right)_{\Omega_h} +\sum_{F\in\mathcal{E}_h^i} \langle (Q^\star \cdot n^+)(Q^+ - Q^-)\cdot \avg{w}\rangle_F \notag                                                                                           \\
                             & - \alpha \Big(\sum_{F\in\mathcal{E}_h^i} h_F^{-1}\langle \jump{Q\cdot n} \jump{w\cdot n}\rangle_F                                      +  \sum_{F\in\mathcal{E}_h^\partial} h_F^{-1}\langle (Q\cdot n) (w\cdot n)\rangle_F \Big)  \label{eqn:f_impl} \\
                             & -\delta_{\text{up}}\sum_{F\in\mathcal{E}_h^i} \langle \left|Q^\star \cdot n^+\right|(Q^+ - Q^-)\cdot (w^+-w^-)\rangle_F\notag                                                                                                                        \\
        f^\expl(w,Q,t)       & := (w\cdot f(t))_{\Omega_h},\label{eqn:f_expl}                                                                                                                                                                                                      
    \end{align}
\end{subequations}
where $\delta_{\text{up}}=1$ for the upwind flux and $\delta_{\text{up}}=0$ for the central flux. The pressure gradient $\widetilde{g}$ and the function $\widetilde{\Gamma}$ that enforces the incompressibility constraint are
\begin{equation}
        \widetilde{g}(w,p)   = \widetilde{\Gamma}(p,w)  := (p \nabla\cdot w)_{\Omega_h} -  \sum_{F\in\mathcal{E}_h^i} \langle \jump{w\cdot n}\avg{p}\rangle_{F}-  \sum_{F\in\mathcal{E}_h^\partial} \langle (w\cdot n)p\rangle_{F} \label{eqn:DG_pressure_gradient}
\end{equation}
In contrast to (2.29) in \cite{Guzman2016}, the boundary condition $n\cdot Q^{n+1}|_{\partial \Omega}=0$ is enforced weakly through the penalty term in the second line of \eqref{eqn:f_impl}, and this introduces additional boundary integrals. However, as the authors of \cite{Guzman2016} stress, enforcing the boundary condition strongly is only necessary for their theoretical analysis, and hence our approach is consistent with their spatial discretisation.  Although in \eqref{eqn:incompressible_euler_momentum} the right-hand side of \eqref{eqn:f_expl} is independent of the velocity, we write down the scheme in a slightly more general form by allowing the $f^\expl$ to also depend on $Q$. This might prove useful in future extensions of our work.
\subsection{Hybridisation}\label{sec:hybridisation}
To hybridise the problem in \eqref{eqn:original_DG_scheme} we introduce an additional function $\lambda\in V_\text{trace}$ defined on the trace space
\begin{equation}
    V_{\text{trace}} = \{ v\in L^2(\mathcal{E}_h) : v|_F \in P_k(F) \;\text{for all facets}\;F\in \mathcal{E}_h \}
\end{equation}
and replace the numerical fluxes in \eqref{eqn:flux_original} by
\begin{equation}
    \begin{aligned}
        \widehat{\sigma}   & := \lambda \Id  + \alpha h_F^{-1} \jump{Q\cdot n} \Id       + \begin{cases}
                                                                                               \avg{Q}\otimes Q^\star       & \text{for the central flux}, \\
                                                                                               Q^{\text{up}}\otimes Q^\star & \text{for the upwind flux}
                                                                                           \end{cases} \\
        \widehat{Q}\cdot n & := Q\cdot n + \tau (p-\lambda).
    \end{aligned}\label{eqn:flux_hdg}
\end{equation}
for some stability parameter $\tau>0$.
The following transmission condition guarantees that the normal component of $\widehat{Q}$ is single-valued:
\begin{equation}
    0 = \langle \jump{\widehat{Q}\cdot n} \mu \rangle_F = \langle (\jump{Q\cdot n} + 2\tau\avg{p-\lambda})\mu \rangle_F \qquad \text{for all $F\in\mathcal{E}_h^i$ and all $\mu\in V_{\text{trace}}$}. \label{eqn:jump_condition_interior}
\end{equation}
The homogeneous Dirichlet boundary condition $n\cdot Q|_{\partial \Omega}=0$ is enforced weakly by requiring that
\begin{equation}
    0 = \langle (\widehat{Q}\cdot n) \mu \rangle_F = \langle (Q\cdot n + \tau(p-\lambda))\mu \rangle_F \qquad \text{for all $F\in\mathcal{E}_h^\partial$ and all $\mu\in V_{\text{trace}}$}.\label{eqn:jump_condition_boundary}
\end{equation}
As a result of the hybridisation, \eqref{eqn:original_DG_scheme} gets replaced by the following problem: given the velocity $Q^n\in V_Q$ at the current timestep find the solution
$(Q^{n+1},p^{n+1},\lambda^{n+1})\in V_Q\times V_p\times V_{\text{trace}}$ at the next timestep such that
\begin{equation}
    \begin{aligned}
        (w\cdot Q^{n+1})_{\Omega_h} - \Delta t \left(f^\impl(w,Q^{n+1},\mathcal{P}(Q^n)) + g(w,p^{n+1},\lambda^{n+1})\right) & = (w\cdot Q^{n})_{\Omega_h} + \Delta t f^\expl(w, Q^n,t^n) \\
        \text{subject to}\qquad\Gamma(\psi,\mu,Q^{n+1},\lambda^{n+1})                                                        & = 0
    \end{aligned}\label{eqn:implicit_HDG_scheme}
\end{equation}
for all test functions $(w,\psi,\mu) \in V_Q\times V_p\times V_{\text{trace}}$. The weak pressure gradient $g$ and the function $\Gamma$ are given by
\begin{subequations}
    \begin{align}
        g(w,p,\lambda)               & := (p \nabla\cdot w)_{\Omega_h} -  \sum_{F\in\mathcal{E}_h^i} \langle \jump{w\cdot n}\lambda\rangle_{F}  -  \sum_{F\in\mathcal{E}_h^\partial} \langle (w\cdot n)\lambda\rangle_{F}                             \\
        \Gamma(\psi,\mu,Q,p,\lambda) & := (\psi \nabla\cdot Q)_{\Omega_h} + 2 \sum_{F\in\mathcal{E}_h^i} \langle \avg{\tau (p-\lambda)\psi}\rangle_{F}            + \sum_{F\in\mathcal{E}_h^\partial} \langle \tau (p-\lambda)\psi\rangle_{F}         \\
                                     & +\sum_{F\in\mathcal{E}_h^i} \langle \left(\jump{Q\cdot n}+2\tau \avg{p-\lambda}\right)\mu \rangle_{F} + \sum_{F\in\mathcal{E}_h^\partial} \langle \left(Q\cdot n+\tau (p-\lambda)\right)\mu \rangle_{F}\notag,
    \end{align}
\end{subequations}
which should be compared to \eqref{eqn:DG_pressure_gradient}.
\subsection{Splitting method}\label{sec:splitting_method}
Finally, instead of solving \eqref{eqn:implicit_HDG_scheme} directly for $Q^{n+1},p^{n+1},\lambda^{n+1}$, we can split this computation into two stages following Chorin's projection method \cite{Chorin1967,Chorin1968}. For this, we first compute a tentative velocity $\overline{Q}$ by dropping the pressure gradient $g$ in the first line of \eqref{eqn:implicit_HDG_scheme} and solving
\begin{equation}
    (w\cdot \overline{Q})_{\Omega_h} - \Delta t f^\impl(w,\overline{Q},\mathcal{P}(Q^n)) = (w\cdot Q^n) + \Delta t f^\expl(w,Q^n)\label{eqn:tentative_velocity}
\end{equation}
for all test functions $w\in V_Q$. The divergence of $\overline{Q}$ is not zero since we do not enforce the incompressibility condition in the second line of \eqref{eqn:implicit_HDG_scheme}. This can be corrected by computing $(\delta Q,p^{n+1},\lambda^{n+1})\in V_Q\times V_p\times V_{\text{trace}}$ which are obtained by solving
\begin{equation}
    \begin{aligned}
        (w\cdot\delta Q)_{\Omega_h} - g(\psi,p^{n+1},\lambda^{n+1})             & = -\frac{1}{\Delta t}\operatorname{Div}(w,\overline{Q}) \\
        \text{subject to}\qquad \Gamma(\psi,\mu,\delta Q,p^{n+1},\lambda^{n+1}) & = 0
    \end{aligned}\label{eqn:split_pressure_solve}
\end{equation}
for all test functions $(w,\psi,\mu) \in V_Q\times V_p\times V_{\text{trace}}$ where the weak divergence is given by
\begin{equation}
    \text{Div}(\psi,Q) := (\psi \nabla \cdot Q)_{\Omega_h}-2\sum_{F\in\mathcal{E}_h^i}\langle(\avg{\psi Q}-\avg{\psi}\avg{Q})\cdot n\rangle_F-\sum_{F\in\mathcal{E}_h^\partial}\langle\psi (Q\cdot n)\rangle_F.\label{eqn:weak_divergence}
\end{equation}
The divergence-free velocity at the next timestep is then given by $Q^{n+1}=\overline{Q}+\Delta t\;\delta Q$. The big advantage of solving \eqref{eqn:split_pressure_solve} instead of \eqref{eqn:implicit_HDG_scheme} is that \eqref{eqn:split_pressure_solve} is the well-known HDG discretisation (see e.g. \cite{Cockburn2009} or \cite[Section 3]{Cockburn2010}) of the mixed Poisson problem
\begin{xalignat}{2}
    \delta Q+\nabla p^{n+1} &= 0, & \nabla\cdot \delta Q &= b \label{eqn:mixed_poisson}
\end{xalignat}
for the right-hand side $b=-\frac{1}{\Delta t}\nabla\cdot \overline{Q}$. Crucially, the linear operator in \eqref{eqn:mixed_poisson} is a special case of the diffusion operator considered in \cite{Cockburn2014}. Because of this, we can use the efficient non-nested multigrid preconditioner proposed in \cite{Cockburn2014} to solve \eqref{eqn:split_pressure_solve}. Obviously, splitting the update $Q^n\mapsto (Q^{n+1},p^{n+1},\lambda^{n+1})$ into two steps will introduce an additional splitting error. Hence, as will be discussed in detail below, the splitting method will not be used on its own. Instead we will employ it to precondition the computation of $(Q^{n+1},p^{n+1},\lambda^{n+1})$ in the full problem \eqref{eqn:implicit_HDG_scheme}.
\subsection{IMEX timesteppers}\label{sec:imex}
So far the above method uses a simple timestepping approach, and we therefore expect it to be first order in time. As described in \cite[Section 4]{Ueckermann2016}, this can be improved by using the higher-order IMEX \cite{Ascher1997,pareschi2005implicit} scheme.

For this, let $(Q,p,\lambda)\in V_Q\times V_p\times V_{\text{trace}}$ be the time dependent state of the system. Consider the time-evolution equation in weak form
\begin{equation}
    (w\cdot \partial_t Q)_{\Omega_h} = f^\impl(w,Q,\mathcal{P}(Q))  + f^\expl(w;t) + g(w,p,\lambda)\label{eqn:time_evolution}
\end{equation}
subject to the weak incompressibility constraint
\begin{equation}
    \Gamma(\psi,\mu,Q,p,\lambda) = 0\label{eqn:constraint}
\end{equation}
where $w\in V_Q$, $\psi\in V_p$, $\mu\in V_{\text{trace}}$ are test functions. Restricting ourselves to the special case $b_0^\impl=a_{i,0}^\impl=0$ as in \cite{Ueckermann2016}, an $s$-stage IMEX-RK scheme is given by: find $(Q^{n+1},\delta p,\delta \lambda)\in V_Q\times V_p\times V_{\text{trace}}$ such that
\begin{equation}
    \begin{aligned}
        (w\cdot Q^{n+1})_{\Omega_h} & = (w\cdot Q^n)_{\Omega_h} + \Delta t \sum_{i=1}^{s-1} b_i^\impl \left(f^\impl(w,Q_i,\mathcal{P}(Q_{i-1}))+g(w,p_i,\lambda_i)\right) \\&\qquad +\Delta t \sum_{i=0}^{s-1} b_i^\expl f^\expl(w;t^n+c_i \Delta t)+\Delta t\;b_{s-1}^\expl\; g(w,\delta p,\delta \lambda)
    \end{aligned}\label{eqn:Q_nplus1_equation}
\end{equation}
subject to the constraint
\begin{equation}
    \Gamma(\psi,\mu,Q^{n+1},\delta p,\delta \lambda) = 0\label{eqn:Q_nplus1_equation_constraint}
\end{equation}
for all test functions $(w,\psi,\mu)\in V_Q\times V_p\times V_{\text{trace}}$. The pressures $\delta p$ and $\delta \lambda$ that are computed when solving \eqref{eqn:final_update} ensure that $Q^{n+1}$ is divergence free. Although the (hybridised) pressure at the next timestep can be computed as $p^{n+1} = p_{s-1}+\delta p$ and $\lambda^{n+1} = \lambda_{s-1}+\delta \lambda$, these quantities will only be correct to first order in the timestep size. To improve on this, a higher-order pressure $p^{n+1}$ (and $\lambda^{n+1}$) at the next timestep can be obtained with an additional linear solve as discussed in Section \ref{sec:pressure_reconstruction}.

The stage variables $Q_1,Q_2,\dots,Q_{s-1}\in V_Q$, $p_1,p_2,\dots,p_{s-1}\in V_p$ and $\lambda_1,\lambda_2,\dots,\lambda_{s-1}\in V_{\text{trace}}$ are obtained as follows: for $i=1,2,\dots,s-1$ find $(Q_i,p_i,\lambda_i)\in V_Q\times V_p\times V_{\text{trace}}$ such that
\begin{equation}
    (Q_i\cdot w)_{\Omega_h} - \Delta t\; a_{i,i}^\impl \big(f^\impl(w,Q_i,\mathcal{P}(Q_{i-1}))+g(w,p_i,\lambda_i)\big)
    = r_i(w)\qquad\text{for all $w\in V_Q$}\label{eqn:imex_implicit_solve}
\end{equation}
with $Q_0=Q^n$, $p_0=p^n$ subject to the constraint
\begin{equation}
    \Gamma(\psi,\mu,Q_i,p_i,\lambda_i) = 0\qquad\text{for all $\psi\in V_p$, $\mu\in V_{\text{trace}}$}.\label{eqn:imex_implicit_solve_constraint}
\end{equation}
At each stage $i$ the residual $r_i$ is a one-form defined by
\begin{equation}
    \begin{aligned}
        r_i(w) & := (w\cdot Q^n)_{\Omega_h} + \Delta t \sum_{j=1}^{i-1} a_{i,j}^\impl \left(f^\impl(w,Q_j,\mathcal{P}(Q_{j-1}))+g(w,p_j,\lambda_j)\right) +\Delta t \sum_{j=0}^{i-1}a_{i,j}^\expl f^\expl(w;t^n+c_j \Delta t).
    \end{aligned}
    \label{eqn:imex_residual_preliminary}
\end{equation}
The coefficients $a_{i,j}^\impl$, $b_i^\impl$, $a_{i,j}^\expl$, $b_i^\expl$ and $c_i$ that define a particular IMEX method are usually written down in the form of Butcher tableaus, see \appref{sec:butcher_tableaus} for details. To avoid re-evaluating $f^\impl$ and to ensure that corresponding terms cancel in the projection method, we follow the strategy in \cite{Ueckermann2016} and use \eqref{eqn:imex_implicit_solve} to express $f^\impl(w,Q_j,\mathcal{P}(Q_{j-1}))$ in terms of $(Q_j\cdot w)_{\Omega_h}$ and $r_j(w)$. With this, we arrive at the following recursive definition of the residual in \eqref{eqn:imex_residual_preliminary}
\begin{equation}
    \begin{aligned}
        r_i(w) & := (w\cdot Q^n)_{\Omega_h} + \sum_{j=1}^{i-1} \frac{a_{i,j}^\impl}{a_{j,j}^\impl} \big((Q_j\cdot w)_{\Omega_h}-r_j(w)\big) +\Delta t \sum_{j=0}^{i-1}a_{i,j}^\expl f^\expl(w;t^n+c_j \Delta t).
    \end{aligned}
    \label{eqn:imex_residual}
\end{equation}
\eqref{eqn:imex_implicit_solve} together with the constraint in \eqref{eqn:imex_implicit_solve_constraint} is a linear equation for $(Q_i,p_i,\lambda_i)$. Instead of solving this linear system in one go as in Section \ref{sec:hybridisation}, we proceed in two stages and use the splitting method from Section \ref{sec:splitting_method} as a preconditioner for a Richardson iteration, which is written down in Alg.~\ref{alg:richardson}. The resulting timestepping algorithm is shown in Alg.~\ref{alg:imex_hdg}. Note that in the calculation of the final residual in \eqref{eqn:final_residual} we again expressed $f^\impl(w,Q_i,\mathcal{P}(Q_{i-1}))$ in terms of $(Q_i\cdot w)_{\Omega_h}$ and $r_i(w)$ (compare the manipulations that transform \eqref{eqn:imex_residual_preliminary} into \eqref{eqn:imex_residual}).
\begin{algorithm}
    \caption{IMEX-HDG timestepping with $s$ stages (optionally based on the projection method): compute $Q^{n+1},p^{n+1}$ at the next timestep given $Q^{n},p^{n}$ at the current time $t^n=n\Delta t$.}
    \label{alg:imex_hdg}
    \begin{algorithmic}[1]
        \State{Initialise $Q_0:=Q^n$, $p_0:=p^n$}
        \For {$i=1,\dots,s-1$}
        \State{Compute the residual $r_i(w)$ defined in \eqref{eqn:imex_residual}}
        \State{Set $Q^\star_{i-1}:=\mathcal{P}(Q_{i-1})$}
        \State{Compute $Q_i,p_i,\lambda_i$ by solving \eqref{eqn:imex_implicit_solve} subject to the constraint in \eqref{eqn:imex_implicit_solve_constraint} with Alg.~\ref{alg:richardson}.}
        \EndFor
        \State{Compute $Q^{n+1}$, $\delta p$, $\delta \lambda$ by solving \eqref{eqn:Q_nplus1_equation} subject to the constraint in \eqref{eqn:Q_nplus1_equation_constraint}. For this, solve
            \begin{equation}
                \begin{aligned}
                    (Q^{n+1}\cdot w)_{\Omega_h} -
                    \Delta t\;b_{s-1}^\impl g(w,\delta p,\delta \lambda) & = r^{n+1}(w) \\
                    \Gamma(\phi,\mu,Q^{n+1},\delta p,\delta \lambda)     & = 0
                \end{aligned}\label{eqn:final_update}
            \end{equation}
        }
        with
        \begin{equation}
            \begin{aligned}
                r^{n+1}(w) & := (w\cdot Q^n)_{\Omega_h} +  \sum_{i=1}^{s-1} \frac{b_i^\impl}{a_{i,i}^\impl} \big((Q_i\cdot w)_{\Omega_h}-r_i(w)\big) +\Delta t \sum_{i=0}^{s-1} b_i^\expl f^\expl(w;t^n+c_i \Delta t) .
            \end{aligned}\label{eqn:final_residual}
        \end{equation}
        \State{Reconstruct the (hybridised) pressure at the next timestep by solving \eqref{eqn:pressure_reconstruction_compact} for $p=p^{n+1}$, $\lambda=\lambda^{n+1}$.}
    \end{algorithmic}
\end{algorithm}
It is worth noting that the simple, first order timestepping method which was used in \cite{Guzman2016} and which is written down for the HDG discretisation Section \ref{sec:hybridisation} can be expressed as a two-stage IMEX method by using the Butcher tableaus in Table \ref{tab:butcher_tableau_imex_euler}.
\begin{algorithm}
    \caption{Preconditioned Richardson iteration. Approximately solve \eqref{eqn:imex_implicit_solve} subject to the constraint in \eqref{eqn:imex_implicit_solve_constraint}, using the projection method as a preconditioner with a fixed number of $n_R$ iterations.}
    \label{alg:richardson}
    \begin{algorithmic}[1]
        \State{Initialise $Q_i^{(0)} = Q_{i-1}$, $p_i^{(0)} = p_{i-1}$, $\lambda_i^{(0)} := \lambda_{i-1}$.}
        \For {$k=1,2,\dots,n_R$}
        \State{Set
            \begin{equation}
                \delta r_i(w) = r_i(w)- (Q_i^{(k-1)}\cdot w)_{\Omega_h} + a_{i,i}^\impl\Delta t \left(f^\impl(w,Q_i^{(k-1)},Q_{i-1}^\star)+g(w,p_i^{(k-1)},\lambda_i^{(k-1)})\right)
            \end{equation}
        }
        \State{\textbf{Step I}: Compute the tentative velocity $\overline{Q}$ by solving (cf. \eqref{eqn:tentative_velocity})
            \begin{equation}
                \begin{aligned}
                    (\overline{Q}\cdot w)_{\Omega_h} - a_{i,i}^\impl\Delta t\; f^\impl(w,\overline{Q},Q_{i-1}^\star) & = \delta r_i(w).
                \end{aligned}
            \end{equation}
        }
        \State{\textbf{Step II}: given $\overline{Q}$, compute increments $\delta p$, $\delta \lambda$ and $\delta Q$ by solving (cf. \eqref{eqn:split_pressure_solve})
            \begin{equation}
                \begin{aligned}
                    (\delta Q\cdot w)_{\Omega_h} - g(w,\delta p,\delta \lambda) & = -\frac{1}{a_{i,i}^\impl \Delta t} \text{Div}(w,\overline{Q}) \\
                    \Gamma(\psi,\mu,\delta Q,\delta p,\delta \lambda)           & = 0
                \end{aligned}\label{eqn:pressure_increment}
            \end{equation}
        }
        \State{\textbf{Step III}: Update
            \begin{xalignat}{3}
                Q_i^{(k)} &= Q_i^{(k-1)} + \overline{Q} + a_{i,i}^\impl \Delta t\;\delta Q, &
                p_i^{(k)} &= p_i^{(k-1)} + \delta p, &
                \lambda_i^{(k)} &= \lambda_i^{(k-1)} + \delta \lambda.
            \end{xalignat}
        }
        \EndFor
        \State{\Return $Q_i^{(n_R)},p_i^{(n_R)},\lambda_i^{(n_R)}$}
    \end{algorithmic}
\end{algorithm}
\subsection{Pressure reconstruction}\label{sec:pressure_reconstruction}
In the final step of Alg.~\ref{alg:imex_hdg} we recover the (hybridised) pressure $p^{n+1}$, $\lambda^{n+1}$ directly from $Q^{n+1}$. This is possible since at any given time the pressure $p$ can be reconstructed from the velocity $Q$ by solving an elliptic problem. To see this, take the divergence of \eqref{eqn:incompressible_euler_momentum} in the domain $\Omega$ and multiply it by the outward normal $n$ on the boundary $\partial \Omega$. Since $\nabla\cdot \partial_t Q = 0$ in $\Omega$ and $n\cdot \partial_t Q=0$ on $\partial \Omega$, this results in the following boundary value problem
\begin{equation}
    \begin{aligned}
        -\Delta p        & = \nabla \cdot f_p\qquad\text{in $\Omega$}    \\
        -n\cdot \nabla p & = n\cdot f_p \qquad\text{on $\partial\Omega$}
    \end{aligned}\label{eqn:pressure_reconstruction_primal}
\end{equation}
with the function
\begin{equation}
    f_p := -f + (Q\cdot \nabla) Q.
\end{equation}
To make the solution $p$ of \eqref{eqn:pressure_reconstruction_primal} unique, we require that $\int_\Omega p\;dx=0$. The expression $n\cdot f_p\vert_{\partial \Omega}$ can be simplified as follows: Since $n\cdot Q\vert_{\partial \Omega}=0$, we can write on the boundary
\begin{equation}
    n\cdot f_p = -n\cdot f + (Q\cdot \nabla) (n\cdot Q)  = -n\cdot f+(Q \cdot \nabla^\parallel)(n\cdot Q) = -n\cdot f
    \label{eqn:n_dot_fp}
\end{equation}
where $\nabla^\parallel=\nabla - n(n\cdot \nabla)$ is the tangential derivative on the boundary. The final identity in \eqref{eqn:n_dot_fp} follows since on the boundary $n\cdot Q=0$ obviously does not change in the tangential direction.

By introducing the new variable $U$, the second order problem in \eqref{eqn:pressure_reconstruction_primal} can be written in mixed form as
\begin{xalignat}{2}
    U + \nabla p &= 0,& \nabla \cdot U &= \nabla \cdot f_p\label{eqn:mixed_poisson_reconstruction}
\end{xalignat}
subject to the boundary condition $n\cdot U\vert_{\partial \Omega} = -n\cdot f$. \eqref{eqn:mixed_poisson_reconstruction} can be solved through hybridisation in exactly the same way as the corresponding problem in \eqref{eqn:mixed_poisson}. The weak formulation of the problem is (cf. \eqref{eqn:split_pressure_solve}): find $(U,p,\lambda) \in V_Q\times V_p\times V_{\text{trace}}$ such that
\begin{equation}
    \begin{aligned}
        (U\cdot w)_{\Omega_h} - g(w,p,\lambda) & = \text{Div}(\psi,f_p)                                            \\
        \Gamma(\psi,\mu,U,p,\lambda)           & = -\sum_{F\in\mathcal{E}_h^\partial}\langle \mu n\cdot f\rangle_F.
    \end{aligned}\label{eqn:pressure_reconstruction_compact}
\end{equation}
Crucially, the bilinear form on the left hand side of \eqref{eqn:pressure_reconstruction_compact} is the same as in \eqref{eqn:final_update} and \eqref{eqn:pressure_increment}. We will discuss an efficient preconditioner for the resulting linear operator in the next section.
\subsection{Solvers and preconditioners}\label{sec:solvers}
To compute the tentative velocity in \eqref{eqn:tentative_velocity} we use a GMRES solver with ILU0 preconditioner. To solve the mixed problem in Eqs. \eqref{eqn:final_update}, \eqref{eqn:pressure_increment} and \eqref{eqn:pressure_reconstruction_compact} we first apply a Schur-complement reduction to reduce the problem to an equation for the facet variables $\lambda$; in Firedrake this is realised with the static condensation preconditioner class \texttt{firedrake.SCPC} \cite{Gibson2020}. The resulting system is then solved with a GMRES iteration, preconditioned with the non-nested multigrid method described in \cite{Cockburn2014} (class \texttt{firedrake.GTMGPC}). On the finest level, a Chebyshev smoother based on the block-diagonal is used such that the matrix which couples the unknowns on each facet is inverted exactly; this is achieved with \texttt{firedrake.ASMStarPC}. For simplicity, the coarse level problem is solved with a standard AMG method with Chebyshev smoothers. Two pre- and post-smoothing iterations are used on all levels. In the non-nested multigrid algorithm the coarse level space consists of piecewise linear conforming elements on the original mesh, i.e. it is formulated in the space
\begin{equation}
    V_{c} = \text{P}_1 := \{ v\in C(\Omega) : v|_K \in P_1(K) \;\text{for all}\;K\in \Omega_h \}.
\end{equation}
The weak form on the coarse level is obtained by re-discretisation of the Laplace operator
\begin{equation}
    a_{c}(\phi_c,\psi_c) = (\nabla \psi_c \cdot \nabla \phi_c)_{\Omega_h}\qquad\text{for all $\phi_c,\psi_c\in V_{c}$}
\end{equation}
and the prolongation $\mathscr{P}:V_c\rightarrow V_{\text{trace}}$ is given by the injection on the skeleton $\mathcal{E}_h$
\begin{xalignat}{2}
    \mathscr{P}&: \phi_c \mapsto \phi &\text{such that}\qquad
    \langle \mu \phi_c \rangle_{\mathcal{E}_h} & =\langle \mu \phi \rangle_{\mathcal{E}_h} \qquad\text{for all $\mu\in V_{\text{trace}}$}.
\end{xalignat}
As will be demonstrated below (see Fig.~\ref{fig:niter}), this non-nested multigrid algorithm is robust with respect to the polynomial degree of the DG discretisation and only shows weak growth in the number of iterations as the resolution increases. The exact PETSc solver options are listed in Tab.~\ref{tab:petsc_solver_options}.
\section{Numerical results}\label{sec:results}
The following numerical experiments demonstrate the accuracy and performance of the HDG-IMEX timestepping methods developed above. Results are shown for two widely used test cases. For the Taylor Green vortex an analytical solution is known, which allows us to demonstrate the convergence of the method for different polynomial degrees and IMEX time-integrators. The second testcase describes a double layer shear flow with a significantly more complicated flow structure; this benchmark is also used in \cite{Guzman2016}.
\subsection{Implementation}
All code was implemented in the Firedrake framework \cite{FiredrakeUserManual}, which supports the efficient composition of solvers and preconditioners \cite{Kirby2017} based on PETSc \cite{petsc-user-ref,petsc-efficient}. This was crucial to construct the sophisticated multigrid solver for the solution of the HDG pressure correction equation as discussed in Section \ref{sec:solvers}. The code is freely available under \url{https://github.com/eikehmueller/IncompressibleEulerHDG} with the particular version that can be used to generate the results in this paper archived as \cite{code_release}.
\subsection{Setup}
Numerical experiments are carried out on two-dimensional domains $\Omega=[0,L]\times[0,L]$ with $L>0$. To construct the mesh $\Omega_h$, the domain $\Omega$ is first subdivided into $n^2$ squares of size $h/\sqrt{2}\times h/\sqrt{2}$ with $h=\sqrt{2}L/n$. Each of these squares is then split into two right-angled triangles with two sides of length $h/\sqrt{2}$ and one side of length $h$. As a result, each grid contains $2n^2$ triangular cells; we will often refer to $n$ as the ``grid size'' in the following. In our numerical experiments we use the upwind flux $\delta_{\text{up}}=1$ and set $\alpha=1$ in \eqref{eqn:f_impl}; the stabilisation parameter in \eqref{eqn:flux_hdg} is kept fixed at $\tau=1$.
\subsubsection{Timestepping methods}
We pick timesteppers such that the temporal order of the method is the same as the polynomial degree $k$ used to discretise the pressure field (recall that according to \eqref{eqn:function_spaces} the components of the velocity field are represented by polynomials of degree $k+1$); in this case it is plausible to expect that the temporal and spatial order of convergence are the same and hence spatial and temporal errors are balanced. The following three setups are considered (the Butcher tableaus of the different IMEX methods are given in \appref{sec:butcher_tableaus}):
\begin{enumerate}
    \item The two-stage implicit-explicit Euler method (Tab.~\ref{tab:butcher_tableau_imex_euler}) with $k=1$
    \item The three-stage SSP2(3,3,2) method (Tab.~\ref{tab:butcher_tableau_ssp2_332}) with $k=2$
    \item The four-stage SSP3(4,3,3) method (Tab.~\ref{tab:butcher_tableau_ssp3_433}) with $k=3$
\end{enumerate}
The number of Richardson iterations in Alg.~\ref{alg:richardson} is set to $n_R=2$ in all numerical experiments.
\subsection{Model systems}
\subsubsection{Taylor Green vortex}\label{sec:taylor_green}
\begin{figure}
    \begin{center}
        \includegraphics[width=0.9\linewidth]{\figdir/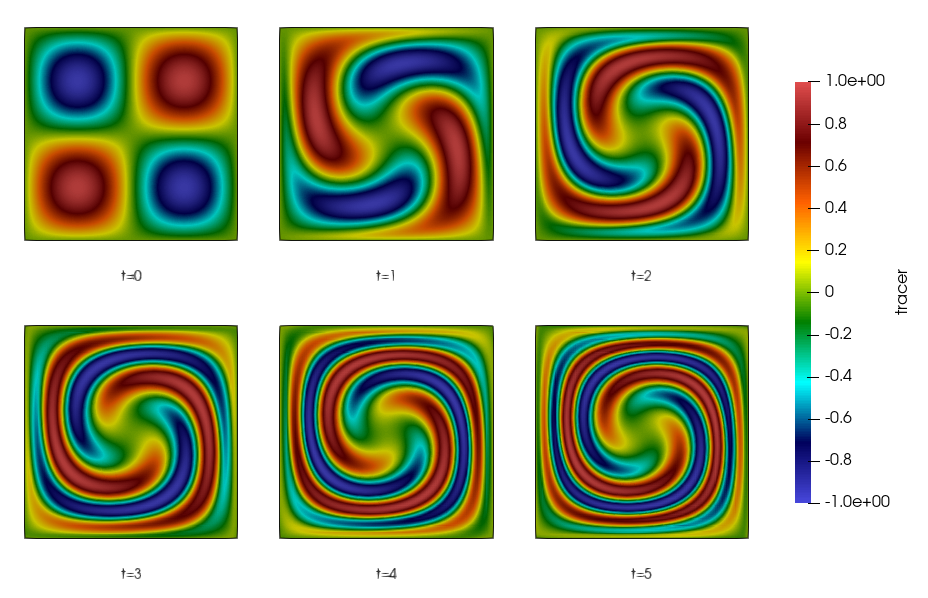}
        \caption{Advection of passive tracer in the Taylor Green vortex at different times on a $32\times 32$ grid with polynomial degree $k=2$. The results were obtained with the IMEX-HDG discretisation and a SSP2(3,3,2) timestepper with timestep size $\Delta t=0.005$. The decay constant in the forcing function in \eqref{eqn:taylor_green_Q0_f} is set to $\kappa=0.1$.}
        \label{fig:tracer_advection}
    \end{center}
\end{figure}

To test the IMEX timestepping methods described in Section \ref{sec:imex} we solve the incompressible Euler equations in \eqref{eqn:incompressible_euler_momentum}, \eqref{eqn:incompressible_euler_constraint} with initial conditions and forcing function chosen such that in the continuum $Q(x,y,t)$ and $p(x,y,t)$ are given by the analytical solution which has the form originally proposed in \cite{Taylor1937}. More specifically, we set
\begin{xalignat}{2}
    Q_0(x,y) &= \begin{pmatrix}-\cos(\frac{2x-1}{2}\pi)\sin(\frac{2y-1}{2}\pi) \\[1ex]
        \sin(\frac{2x-1}{2}\pi)\cos(\frac{2y-1}{2}\pi)\end{pmatrix},&
    f(x,y,t) &= -\kappa e^{-\kappa t} Q_0(x,y)\label{eqn:taylor_green_Q0_f}
\end{xalignat}
in the domain $\Omega=[0,1]\times[0,1]$ which results in the exact solution
\begin{xalignat}{2}
    Q(x,y,t) &= e^{-\kappa t} Q_0(x,y), &
    p(x,y,t) &= e^{-2\kappa t}\left(\frac{4}{\pi^2} - \cos(\frac{2x-1}{2}\pi)\cos(\frac{2y-1}{2}\pi)\right).\label{eqn:exact_solution}
\end{xalignat}
Fig.~\ref{fig:tracer_advection} shows the advection of a passive tracer in the simulated velocity field which was generated with the HDG-IMEX method described in Section \ref{sec:methods}. The results were obtained with a SSP2(3,3,2) timestepper, using polynomial degree $k=2$ for both the pressure and passive tracer, the velocity field has degree $k+1=3$. Tracer advection is treated purely explicitly with a DG-upwinding scheme in the same IMEX timestepper, see \appref{sec:DG_tracer_advection} for further details. Fig.~\ref{fig:tracer_advection} should be compared to the corresponding Fig.~\ref{fig:tracer_advection_dg} in the appendix, which was obtained with the original, fully implicit DG scheme in \cite{Guzman2016} as described in Section \ref{sec:dg_discretisation}.
To ensure that the explicit tracer advection scheme is stable, it was necessary to use a relatively small timestep size of $\Delta t=0.005$ to produce Figs.~\ref{fig:tracer_advection} and \ref{fig:tracer_advection_dg}. However, as the results in Section \ref{sec:convergence} show, the computation of the velocity and pressure field can be performed with an approximately six times larger timestep (see Tab.~\ref{tab:stepsizes}). While it would be possible to avoid this constraint on $\Delta t$ by using a more sophistiacted numerical scheme for the tracer advection or by treating both processes on different timescales, this is not the focus of the present paper.
\subsubsection{Double-layer shear flow}
\begin{figure}
    \begin{center}
        \includegraphics[width=0.9\linewidth]{\figdir/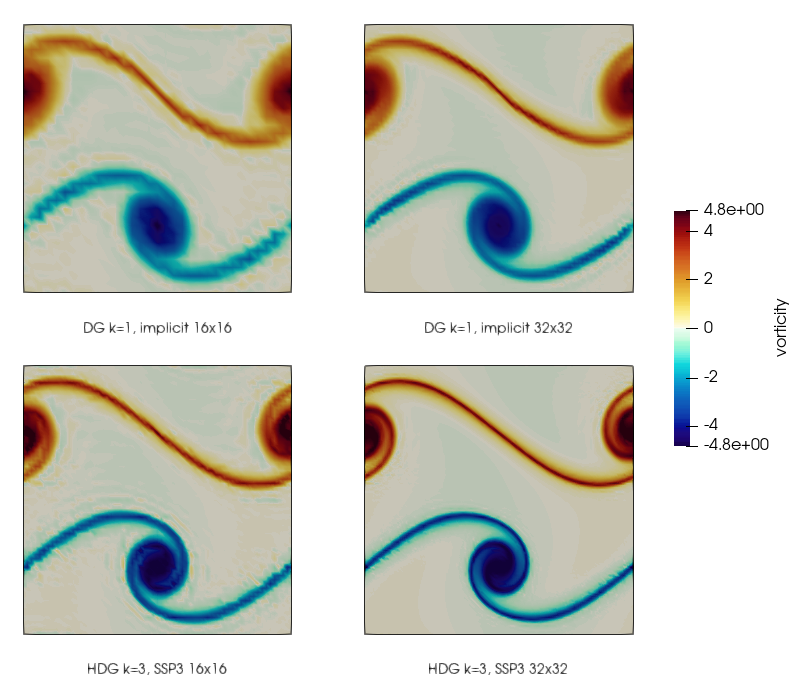}
        \caption{Vorticity for double layer shear flow at time $t=6$ on two different grids. The two upper plots show results obtained with the original method in \cite{Guzman2016}.}
        \label{fig:shear_flow_t6.0}
    \end{center}
\end{figure}

\begin{figure}
    \begin{center}
        \includegraphics[width=0.9\linewidth]{\figdir/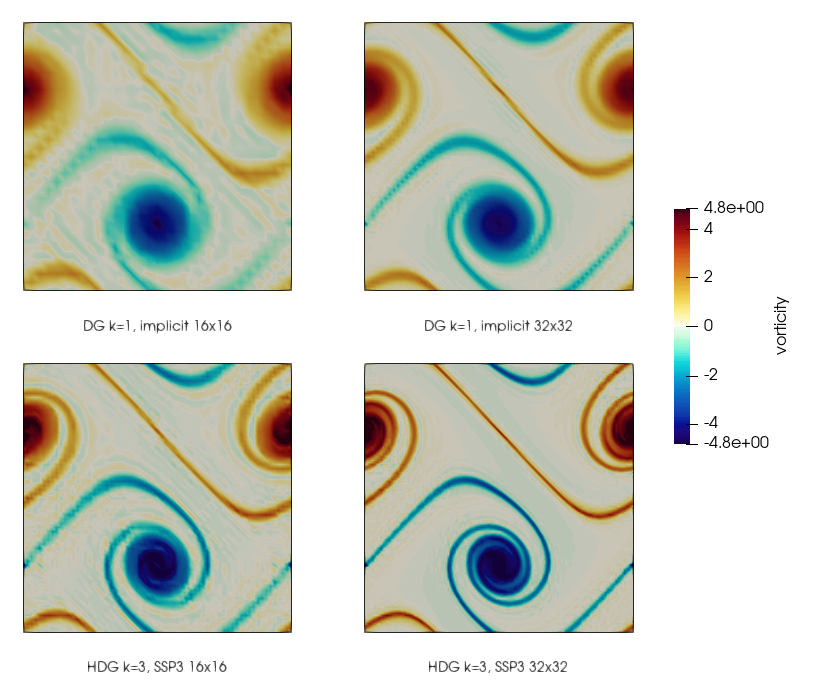}
        \caption{Vorticity for double layer shear flow at time $t=8$ on two different grids. The two upper plots show results obtained with the original method in \cite{Guzman2016}.}
        \label{fig:shear_flow_t8.0}
    \end{center}
\end{figure}

We also consider the double layer shear flow testcase \cite{Bell1989} in the same setup as in \cite{Guzman2016}. For this, the problem is solved in the domain $\Omega=[0,2\pi]\times[0,2\pi]$ and the initial condition is set to
\begin{equation}
    Q_0(x,y)   = \begin{pmatrix}
        Q_{0,x}(y) \\
        \delta \sin(x)
    \end{pmatrix}    \qquad\text{with\quad $Q_{0,x}(y) = \begin{cases}
                \tanh\left(\left(y-\tfrac{\pi}{2}\right)/\rho\right)  & \text{for $0\le y\le \pi$,} \\
                \tanh\left(\left(\tfrac{3\pi}{2}-y\right)/\rho\right) & \text{for $\pi<y\le 2\pi$.}
            \end{cases}$}
    \label{eqn:shear_flow_velocity_IC}
\end{equation}
As in \cite{Guzman2016} the parameters are fixed to $\rho=\pi/15$ and $\delta=0.05$. Periodic boundary conditions are assumed in both dimensions. Figs.~\ref{fig:shear_flow_t6.0} and \ref{fig:shear_flow_t8.0} visualise the vorticity at times $t=6$ and $t=8$ at two different spatial resolutions; these two figures should be compared to Fig. 8 and Fig. 11 in \cite{Guzman2016} respectively. Results are shown for both the original DG method and degree $k=1$ and for the HDG method with a SSP3(4,3,3) timestepper and $k=3$. The upwind flux and a timestep size of $\Delta t=0.04$ was used in all cases.
\subsection{Convergence}\label{sec:convergence}
\begin{figure}
    \begin{center}
        \includegraphics[width=1.0\linewidth]{\figdir/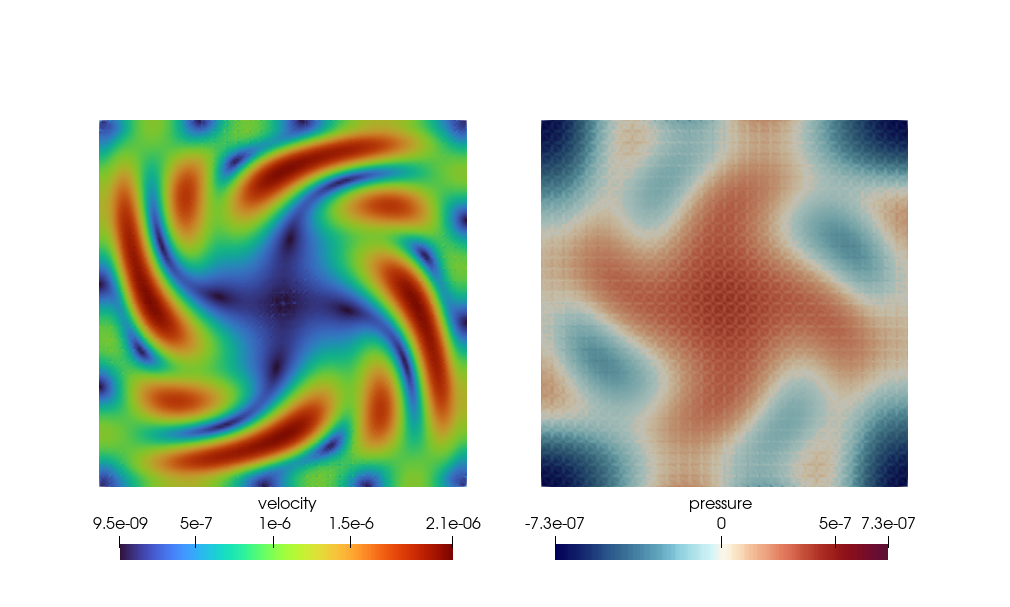}
        \caption{Spatial distribution of the error in velocity (left) and pressure (right) on a $32\times 32$ grid. Results are obtained with a SSP3(4,3,3) timestepper and a polynomial degree of $k=3$.}
        \label{fig:spatial_error}
    \end{center}
\end{figure}

\begin{figure}
    \begin{center}
        \includegraphics[width=1.0\linewidth]{\figdir/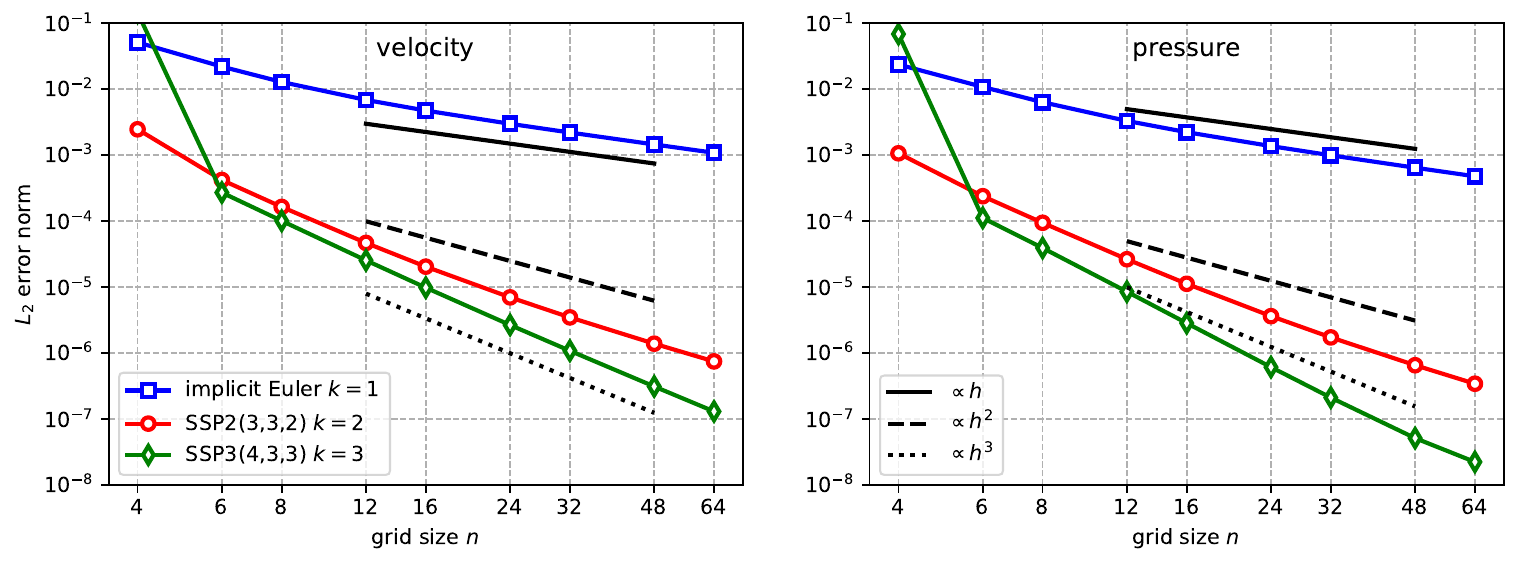}
        \caption{$L_2$ error norm in velocity (left) and pressure (right) for different grid sizes and polynomial degrees.}
        \label{fig:L2error}
    \end{center}
\end{figure}

We verify that the numerical solution of the Taylor Green problem in Section~\ref{sec:taylor_green} converges to the exact solution given in Eqs. \eqref{eqn:exact_solution}. We set $\kappa=0.5$ in all numerical experiments in this section and fix the final time to $T=1$. The number of timesteps is set to $n_t=n$ which results in a timestep size of $\Delta t=T/n_t = L/n = h/\sqrt{2}$; numerical values for the parameters on different grids are given in Tab.~\ref{tab:stepsizes}. As the final six columns of this table shows, the number of trace-unknowns (= size of $V_\text{trace}$) is significantly smaller than the number of original DG unknowns (= size of $V_Q\times V_p$), in particular for higher discretisation orders $k$.
\begin{table}
    \begin{center}
        
    \begin{tabular}{|c|x{1.5cm}|x{1.5cm}|x{1.5cm}|x{1.5cm}|x{1.5cm}|x{1.5cm}|x{1.5cm}|x{1.5cm}|}
    \hline
    \multirow{2}{*}{grid} & \multirow{2}{*}{$h$} & \multirow{2}{*}{$\Delta t$} & \multicolumn{3}{c|}{number of DG unknowns} & \multicolumn{3}{c|}{number of trace unknowns} \\
    &  & & $k=1$ & $k=2$ & $k=3$ & $k=1$ & $k=2$ & $k=3$    \\
    \hline\hline
    
$4\times4$ & $0.353553$ & $0.250000$ & $480$& $832$& $1280$& $112$& $168$& $224$\\\hline
$6\times6$ & $0.235702$ & $0.166667$ & $1080$& $1872$& $2880$& $240$& $360$& $480$\\\hline
$8\times8$ & $0.176777$ & $0.125000$ & $1920$& $3328$& $5120$& $416$& $624$& $832$\\\hline
$12\times12$ & $0.117851$ & $0.083333$ & $4320$& $7488$& $11520$& $912$& $1368$& $1824$\\\hline
$16\times16$ & $0.088388$ & $0.062500$ & $7680$& $13312$& $20480$& $1600$& $2400$& $3200$\\\hline
$24\times24$ & $0.058926$ & $0.041667$ & $17280$& $29952$& $46080$& $3552$& $5328$& $7104$\\\hline
$32\times32$ & $0.044194$ & $0.031250$ & $30720$& $53248$& $81920$& $6272$& $9408$& $12544$\\\hline
$48\times48$ & $0.029463$ & $0.020833$ & $69120$& $119808$& $184320$& $14016$& $21024$& $28032$\\\hline
$64\times64$ & $0.022097$ & $0.015625$ & $122880$& $212992$& $327680$& $24832$& $37248$& $49664$\\\hline
\end{tabular}

        \caption{Grids $n\times n$ with corresponding grid spacing $h$ and timestep size $\Delta t$ (both truncated to six decimal places) used in the numerical experiments. Columns 4-6 show the total number of velocity- and pressure unknowns in the DG discretisation \eqref{eqn:function_spaces} for different polynomial degrees $k$. The corresponding number of trace unknowns is tabulated in the final three columns.}
        \label{tab:stepsizes}
    \end{center}
\end{table}
To illustrate that errors not localised in parts of the domain, Fig.~\ref{fig:spatial_error} shows the spatial distribution of the error on a $32\times 32$ grid with a SSP3(4,3,3) timestepper and $k=3$. The quantity that is plotted is $\|Q(x,T)-Q_{\text{exact}}(x,T)\|_2$ for the velocity and $|p(x,T)-p_{\text{exact}}(x,T)|$ for the pressure.
To quantify convergence with the grid spacing $h$, we compute the $L_2$ error norm
which is
\begin{equation}
    \begin{aligned}
    \|p-p_\text{exact}\|_2 &= \left(\int_\Omega (p(x,T)-p_\text{exact}(x,T))^2\;dx \right)^{1/2},\\
    \|Q-Q_\text{exact}\|_2 &= \left(\int_\Omega \|Q(x,T)-Q_\text{exact}(x,T)\|_2^2\;dx \right)^{1/2}
\end{aligned} 
\end{equation}
for pressure and velocity respectively. Fig.~\ref{fig:L2error} shows the dependence of the error in the $L_2$ norm on the grid size for different polynomial degrees $k$ as the resolution increases. For both velocity and pressure we empirically observe convergence with the asymptotic rate  $\propto h^{k}\propto (\Delta t)^k$.
\subsection{Performance}
Next, we investigate the performance of the different time integrators. For this, we first identify the computational bottlenecks. In the following we will set the final time to $T=1$ for the Taylor Green vortex and to $T=8$ for the double layer shear flow testcase. Again, the number of timesteps is set to $n_t=n$ to balance spatial and temporal errors, which results in the timestep sizes $\Delta t$ shown in Tab.~\ref{tab:stepsizes}.
\subsubsection{Breakdown of time per iteration}
\begin{figure}
    \begin{center}
        \includegraphics[width=1.0\linewidth]{\figdir/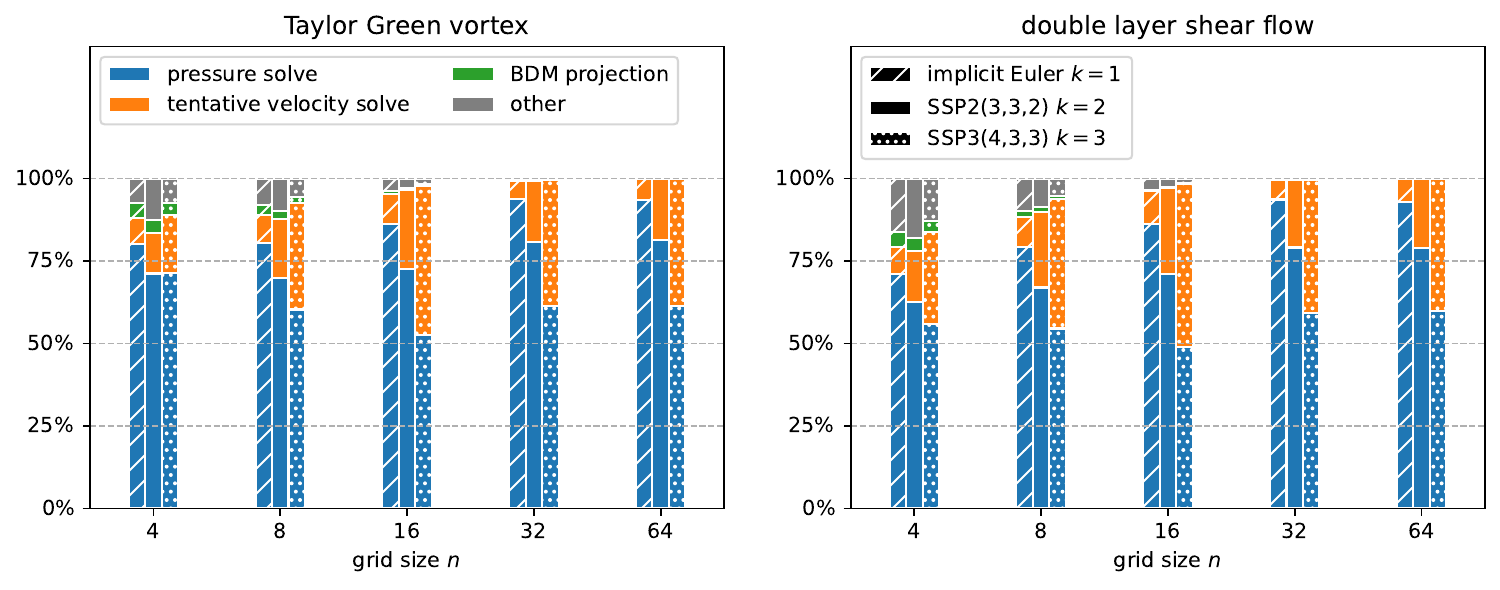}
        \caption{Fraction of runtime spent in the different solver components for varying grid size and different polynomial degrees. Results are shown both for the Taylor Green vortex (left) and the double layer shear flow(right).}
        \label{fig:titer_breakdown}
    \end{center}
\end{figure}

Recall that an $s$-stage method with $n_R$ Richardson iterations requires $s-1$ projections of the velocity to the $\text{BDM}_{k+1}^0$ space, $n_R(s-1)$ solves for the tentative velocity and $n_R(s-1)+2$ mixed pressure solves. We expect these operations to dominate the total cost and this is confirmed by Fig.~\ref{fig:titer_breakdown} which shows a breakdown of the runtime for both testcases; more detailed results can be found in Tabs.~\ref{tab:time_per_iteration_breakdown_taylorgreen} and \ref{tab:time_per_iteration_breakdown_shear} where measurements are shown for the 2-stage implicit Euler method, the 3-stage SSP2(3,3,2) integrator and the 4-stage SSP3(4,3,3) timestepper for both testcases.
\begin{table}
\begin{center}
implicit Euler, $k=1$\\[1ex]
\begin{tabular}{|c|rcrcr|rcrcr|rcrcr|r|}
\hline
grid                                           & 
\multicolumn{5}{|c|}{pressure solve}           &
\multicolumn{5}{|c|}{tentative velocity solve} &
\multicolumn{5}{|c|}{BDM projection}           &
timestep\\
\hline\hline
$ 4\times 4$ &  4 & $\times$ &  0.06 & $=$ &  0.24 &  2 & $\times$ &  0.01 & $=$ &  0.02 &  1 & $\times$ & 0.013 & $=$ & 0.013 &  0.30 \\
$ 6\times 6$ &  4 & $\times$ &  0.05 & $=$ &  0.19 &  2 & $\times$ &  0.01 & $=$ &  0.02 &  1 & $\times$ & 0.009 & $=$ & 0.009 &  0.24 \\
$ 8\times 8$ &  4 & $\times$ &  0.05 & $=$ &  0.20 &  2 & $\times$ &  0.01 & $=$ &  0.02 &  1 & $\times$ & 0.007 & $=$ & 0.007 &  0.24 \\
$12\times12$ &  4 & $\times$ &  0.07 & $=$ &  0.28 &  2 & $\times$ &  0.01 & $=$ &  0.03 &  1 & $\times$ & 0.005 & $=$ & 0.005 &  0.33 \\
$16\times16$ &  4 & $\times$ &  0.11 & $=$ &  0.43 &  2 & $\times$ &  0.02 & $=$ &  0.04 &  1 & $\times$ & 0.005 & $=$ & 0.005 &  0.50 \\
$24\times24$ &  4 & $\times$ &  0.35 & $=$ &  1.39 &  2 & $\times$ &  0.05 & $=$ &  0.09 &  1 & $\times$ & 0.004 & $=$ & 0.004 &  1.50 \\
$32\times32$ &  4 & $\times$ &  0.73 & $=$ &  2.90 &  2 & $\times$ &  0.08 & $=$ &  0.17 &  1 & $\times$ & 0.003 & $=$ & 0.003 &  3.10 \\
$48\times48$ &  4 & $\times$ &  1.64 & $=$ &  6.56 &  2 & $\times$ &  0.21 & $=$ &  0.42 &  1 & $\times$ & 0.003 & $=$ & 0.003 &  7.00 \\
$64\times64$ &  4 & $\times$ &  2.98 & $=$ & 11.92 &  2 & $\times$ &  0.41 & $=$ &  0.82 &  1 & $\times$ & 0.004 & $=$ & 0.004 & 12.76 \\
\hline\end{tabular}
\end{center}

\begin{center}
SSP2(3,3,2), $k=2$\\[1ex]
\begin{tabular}{|c|rcrcr|rcrcr|rcrcr|r|}
\hline
grid                                           & 
\multicolumn{5}{|c|}{pressure solve}           &
\multicolumn{5}{|c|}{tentative velocity solve} &
\multicolumn{5}{|c|}{BDM projection}           &
timestep\\
\hline\hline
$ 4\times 4$ &  6 & $\times$ &  0.05 & $=$ &  0.28 &  4 & $\times$ &  0.01 & $=$ &  0.05 &  2 & $\times$ & 0.008 & $=$ & 0.015 &  0.40 \\
$ 6\times 6$ &  6 & $\times$ &  0.04 & $=$ &  0.24 &  4 & $\times$ &  0.01 & $=$ &  0.05 &  2 & $\times$ & 0.006 & $=$ & 0.012 &  0.35 \\
$ 8\times 8$ &  6 & $\times$ &  0.05 & $=$ &  0.28 &  4 & $\times$ &  0.02 & $=$ &  0.07 &  2 & $\times$ & 0.005 & $=$ & 0.010 &  0.40 \\
$12\times12$ &  6 & $\times$ &  0.08 & $=$ &  0.45 &  4 & $\times$ &  0.03 & $=$ &  0.14 &  2 & $\times$ & 0.004 & $=$ & 0.008 &  0.63 \\
$16\times16$ &  6 & $\times$ &  0.12 & $=$ &  0.72 &  4 & $\times$ &  0.06 & $=$ &  0.24 &  2 & $\times$ & 0.003 & $=$ & 0.007 &  1.00 \\
$24\times24$ &  6 & $\times$ &  0.39 & $=$ &  2.36 &  4 & $\times$ &  0.14 & $=$ &  0.58 &  2 & $\times$ & 0.003 & $=$ & 0.006 &  2.97 \\
$32\times32$ &  6 & $\times$ &  0.82 & $=$ &  4.94 &  4 & $\times$ &  0.28 & $=$ &  1.14 &  2 & $\times$ & 0.003 & $=$ & 0.006 &  6.11 \\
$48\times48$ &  6 & $\times$ &  1.84 & $=$ & 11.05 &  4 & $\times$ &  0.68 & $=$ &  2.71 &  2 & $\times$ & 0.003 & $=$ & 0.007 & 13.80 \\
$64\times64$ &  6 & $\times$ &  3.39 & $=$ & 20.37 &  4 & $\times$ &  1.15 & $=$ &  4.61 &  2 & $\times$ & 0.004 & $=$ & 0.008 & 25.03 \\
\hline\end{tabular}
\end{center}

\begin{center}
SSP3(4,3,3), $k=3$\\[1ex]
\begin{tabular}{|c|rcrcr|rcrcr|rcrcr|r|}
\hline
grid                                           & 
\multicolumn{5}{|c|}{pressure solve}           &
\multicolumn{5}{|c|}{tentative velocity solve} &
\multicolumn{5}{|c|}{BDM projection}           &
timestep\\
\hline\hline
$ 4\times 4$ &  8 & $\times$ &  0.05 & $=$ &  0.39 &  6 & $\times$ &  0.02 & $=$ &  0.10 &  3 & $\times$ & 0.007 & $=$ & 0.020 &  0.55 \\
$ 6\times 6$ &  8 & $\times$ &  0.04 & $=$ &  0.34 &  6 & $\times$ &  0.02 & $=$ &  0.14 &  3 & $\times$ & 0.005 & $=$ & 0.014 &  0.54 \\
$ 8\times 8$ &  8 & $\times$ &  0.05 & $=$ &  0.43 &  6 & $\times$ &  0.04 & $=$ &  0.23 &  3 & $\times$ & 0.004 & $=$ & 0.013 &  0.71 \\
$12\times12$ &  8 & $\times$ &  0.09 & $=$ &  0.72 &  6 & $\times$ &  0.09 & $=$ &  0.54 &  3 & $\times$ & 0.003 & $=$ & 0.010 &  1.31 \\
$16\times16$ &  8 & $\times$ &  0.15 & $=$ &  1.18 &  6 & $\times$ &  0.17 & $=$ &  1.02 &  3 & $\times$ & 0.003 & $=$ & 0.010 &  2.24 \\
$24\times24$ &  8 & $\times$ &  0.44 & $=$ &  3.52 &  6 & $\times$ &  0.41 & $=$ &  2.45 &  3 & $\times$ & 0.003 & $=$ & 0.009 &  6.03 \\
$32\times32$ &  8 & $\times$ &  0.91 & $=$ &  7.27 &  6 & $\times$ &  0.75 & $=$ &  4.52 &  3 & $\times$ & 0.003 & $=$ & 0.010 & 11.86 \\
$48\times48$ &  8 & $\times$ &  2.04 & $=$ & 16.32 &  6 & $\times$ &  1.74 & $=$ & 10.43 &  3 & $\times$ & 0.005 & $=$ & 0.015 & 26.83 \\
$64\times64$ &  8 & $\times$ &  3.77 & $=$ & 30.18 &  6 & $\times$ &  3.14 & $=$ & 18.84 &  3 & $\times$ & 0.007 & $=$ & 0.020 & 49.13 \\
\hline\end{tabular}
\end{center}

    \begin{center}
        \caption{Breakdown of the time spent in the solver components for each timestep for the Taylor Green vortex testcase. For each component, the time per call and the number of calls per timestep is given. The total time per timestep is listed in the final column.}
        \label{tab:time_per_iteration_breakdown_taylorgreen}
    \end{center}
\end{table}
\begin{table}
\begin{center}
implicit Euler, $k=1$\\[1ex]
\begin{tabular}{|c|rcrcr|rcrcr|rcrcr|r|}
\hline
grid                                           & 
\multicolumn{5}{|c|}{pressure solve}           &
\multicolumn{5}{|c|}{tentative velocity solve} &
\multicolumn{5}{|c|}{BDM projection}           &
timestep\\
\hline\hline
$ 4\times 4$ &  4 & $\times$ &  0.01 & $=$ &  0.05 &  2 & $\times$ &  0.00 & $=$ &  0.01 &  1 & $\times$ & 0.003 & $=$ & 0.003 &  0.07 \\
$ 8\times 8$ &  4 & $\times$ &  0.02 & $=$ &  0.09 &  2 & $\times$ &  0.01 & $=$ &  0.01 &  1 & $\times$ & 0.002 & $=$ & 0.002 &  0.12 \\
$16\times16$ &  4 & $\times$ &  0.08 & $=$ &  0.34 &  2 & $\times$ &  0.02 & $=$ &  0.04 &  1 & $\times$ & 0.002 & $=$ & 0.002 &  0.39 \\
$32\times32$ &  4 & $\times$ &  0.67 & $=$ &  2.67 &  2 & $\times$ &  0.08 & $=$ &  0.17 &  1 & $\times$ & 0.002 & $=$ & 0.002 &  2.85 \\
$64\times64$ &  4 & $\times$ &  2.90 & $=$ & 11.62 &  2 & $\times$ &  0.43 & $=$ &  0.86 &  1 & $\times$ & 0.003 & $=$ & 0.003 & 12.50 \\
\hline\end{tabular}
\end{center}

\begin{center}
SSP2(3,3,2), $k=2$\\[1ex]
\begin{tabular}{|c|rcrcr|rcrcr|rcrcr|r|}
\hline
grid                                           & 
\multicolumn{5}{|c|}{pressure solve}           &
\multicolumn{5}{|c|}{tentative velocity solve} &
\multicolumn{5}{|c|}{BDM projection}           &
timestep\\
\hline\hline
$ 4\times 4$ &  6 & $\times$ &  0.01 & $=$ &  0.07 &  4 & $\times$ &  0.00 & $=$ &  0.02 &  2 & $\times$ & 0.002 & $=$ & 0.005 &  0.12 \\
$ 8\times 8$ &  6 & $\times$ &  0.03 & $=$ &  0.16 &  4 & $\times$ &  0.01 & $=$ &  0.06 &  2 & $\times$ & 0.002 & $=$ & 0.004 &  0.24 \\
$16\times16$ &  6 & $\times$ &  0.10 & $=$ &  0.62 &  4 & $\times$ &  0.06 & $=$ &  0.23 &  2 & $\times$ & 0.002 & $=$ & 0.004 &  0.88 \\
$32\times32$ &  6 & $\times$ &  0.76 & $=$ &  4.57 &  4 & $\times$ &  0.29 & $=$ &  1.18 &  2 & $\times$ & 0.002 & $=$ & 0.004 &  5.78 \\
$64\times64$ &  6 & $\times$ &  3.28 & $=$ & 19.66 &  4 & $\times$ &  1.31 & $=$ &  5.24 &  2 & $\times$ & 0.004 & $=$ & 0.008 & 24.94 \\
\hline\end{tabular}
\end{center}

\begin{center}
SSP3(4,3,3), $k=3$\\[1ex]
\begin{tabular}{|c|rcrcr|rcrcr|rcrcr|r|}
\hline
grid                                           & 
\multicolumn{5}{|c|}{pressure solve}           &
\multicolumn{5}{|c|}{tentative velocity solve} &
\multicolumn{5}{|c|}{BDM projection}           &
timestep\\
\hline\hline
$ 4\times 4$ &  8 & $\times$ &  0.01 & $=$ &  0.11 &  6 & $\times$ &  0.01 & $=$ &  0.06 &  3 & $\times$ & 0.002 & $=$ & 0.006 &  0.20 \\
$ 8\times 8$ &  8 & $\times$ &  0.04 & $=$ &  0.29 &  6 & $\times$ &  0.04 & $=$ &  0.21 &  3 & $\times$ & 0.002 & $=$ & 0.006 &  0.54 \\
$16\times16$ &  8 & $\times$ &  0.13 & $=$ &  1.06 &  6 & $\times$ &  0.18 & $=$ &  1.08 &  3 & $\times$ & 0.002 & $=$ & 0.006 &  2.17 \\
$32\times32$ &  8 & $\times$ &  0.88 & $=$ &  7.02 &  6 & $\times$ &  0.80 & $=$ &  4.78 &  3 & $\times$ & 0.003 & $=$ & 0.009 & 11.85 \\
$64\times64$ &  8 & $\times$ &  3.91 & $=$ & 31.27 &  6 & $\times$ &  3.47 & $=$ & 20.79 &  3 & $\times$ & 0.007 & $=$ & 0.020 & 52.16 \\
\hline\end{tabular}
\end{center}

    \begin{center}
        \caption{Breakdown of the time spent in the solver components for each timestep for the double layer shear flow testcase. For each component, the time per call and the number of calls per timestep is given. The total time per timestep is listed in the final column. }
        \label{tab:time_per_iteration_breakdown_shear}
    \end{center}
\end{table}
The results show that the runtime is dominated by the solution of the mixed pressure system. The tentative velocity solve takes up an increasingly larger proportion of the remaining time for higher-order discretisations. For higher resolutions, the results are comparable between the two testcases, which indicates that the more sophisticated structure of the flow in the double layer shear flow does not adversely affect the runtime.
\subsubsection{Robustness of the solvers}\label{sec:solver_robustness}
\begin{figure}
    \begin{center}
        \includegraphics[width=1.0\linewidth]{\figdir/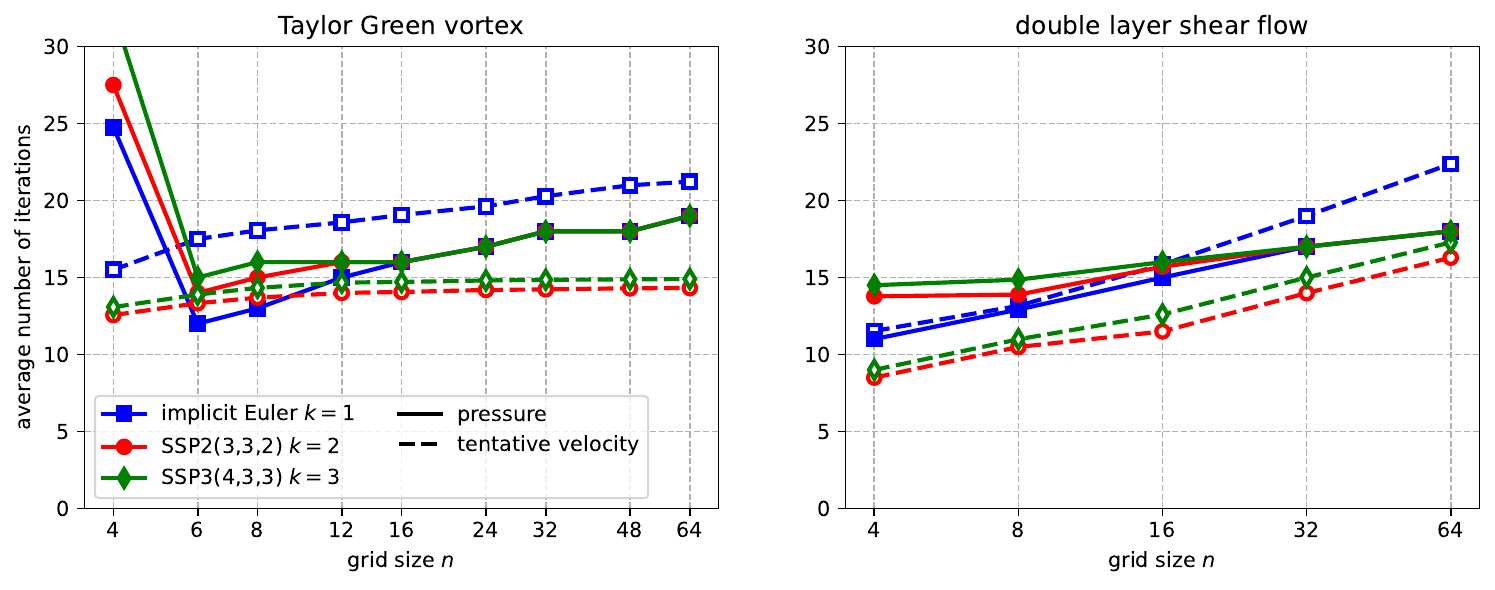}
        \caption{Number of solver as a function of the grid size for different polynomial degrees. Results are shown both for the Taylor Green vortex (left) and the double layer shear flow(right).}
        \label{fig:niter}
    \end{center}
\end{figure}

To interpret the results in the previous section, we investigate the $h$- and $p$- robustness of the pressure- and tentative velocity solves. For this we study the dependence of the number of iterative solver iterations on the grid resolution and the polynomial degree $k$. As in all previous experiments, we manually set relative tolerances of $\varepsilon=10^{-12}$ on the preconditioned residual for the pressure solve and $\varepsilon=10^{-10}$ for the tentative velocity solver. Fig.~\ref{fig:niter} shows the average number of iterations as a function of the grid size both for the tentative velocity solver and the pressure solver; results are shown for both testcases. For the pressure solver, we only consider data collected for the final pressure solve in each timestep. However, we verified that the other solves require approximately the same number of iterations. The figure demonstrates that the performance of the pressure solver is independent of the polynomial degree $k$ as the resolution increases. However, the number of pressure solver iterations shows a modest increase as the resolution is refined. Interestingly, the number of iterations for the tentative velocity solve decreases for higher polynomial degrees $k$. For this solver the dependence of the convergence rate on the grid resolution is different for the two test cases: while for the Taylor Green vortex the number of iterations remains essentially fixed, there is some modest increase for the double layer shear flow benchmark. This is likely related to the more complicated structure of the flow in this case, which shows detail at small length scales (see Figs.~\ref{fig:shear_flow_t6.0} and \ref{fig:shear_flow_t8.0}). In an exploratory study (not shown here) we have also replaced the ILU preconditioner by other methods such as AMG, but this has always led to poorer overall performance.
\subsubsection{Cost per time step}\label{sec:cost_per_timestep}
\begin{figure}
    \begin{center}
        \includegraphics[width=1.0\linewidth]{\figdir/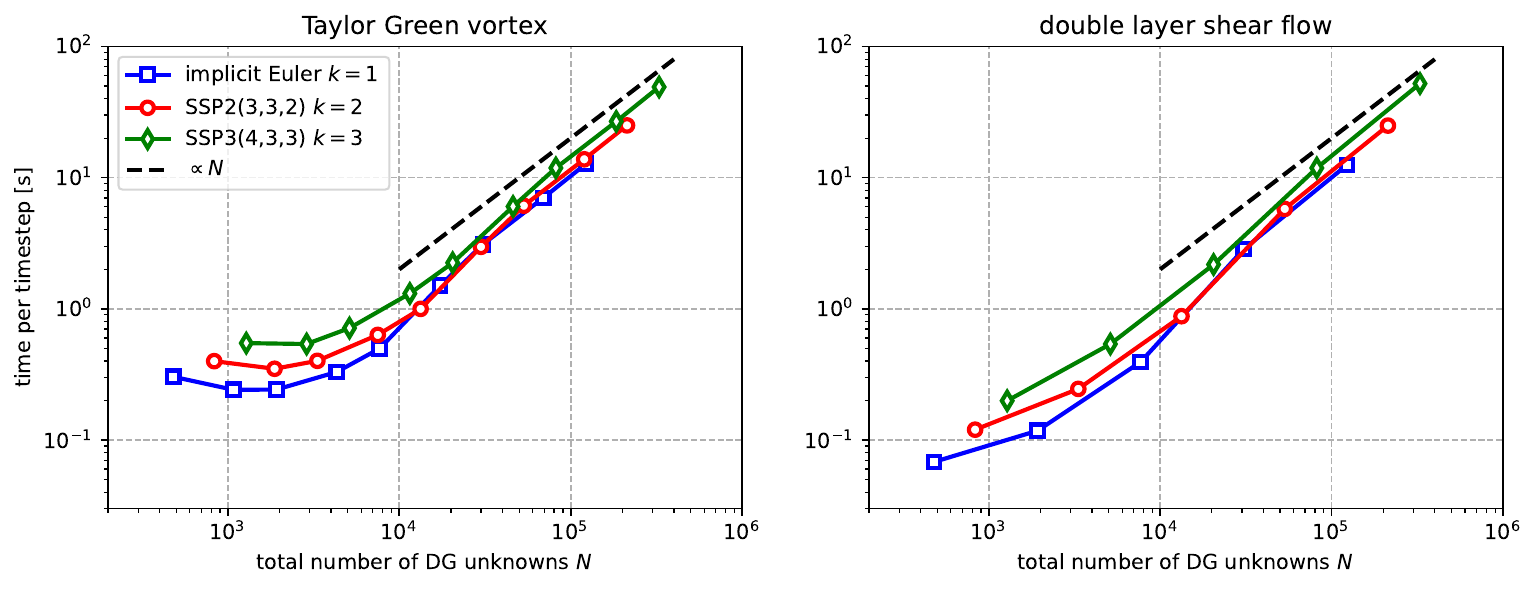}
        \caption{Time per timestep as a function of the total number of DG unknowns (see Tab.~\ref{tab:stepsizes}) for different polynomial degrees. Results are shown both for the Taylor Green vortex (left) and the double layer shear flow(right).}
        \label{fig:titer}
    \end{center}
\end{figure}

Since (ideally) each solver iteration incurs a cost proportional to the number of unknowns, we expect the cost of a single timestep to be proportional to number of unknowns for a given timestepper and polynomial degree $k$.
This time per timestep is plotted as a function of the number of DG unknowns $N$ (cf. Tab.~\ref{tab:stepsizes}) for different polynomial degrees in Fig.~\ref{fig:titer}. Asymptotically the time per timestep grows approximately linearly in $N$. This is consistent with the results Section \ref{sec:solver_robustness}, from which we would expect a slightly faster than linear growth due to the small increase in the number of solver iterations at higher resolution. It should be stressed that this linear growth in runtime that we achieved here is a non-trivial result: naively implicit methods suffer from a rapid growth in cost since the time spent in linear solver increases significantly at higher resolutions.
\subsubsection{Total runtime}
\begin{figure}
    \begin{center}
        \includegraphics[width=1.0\linewidth]{\figdir/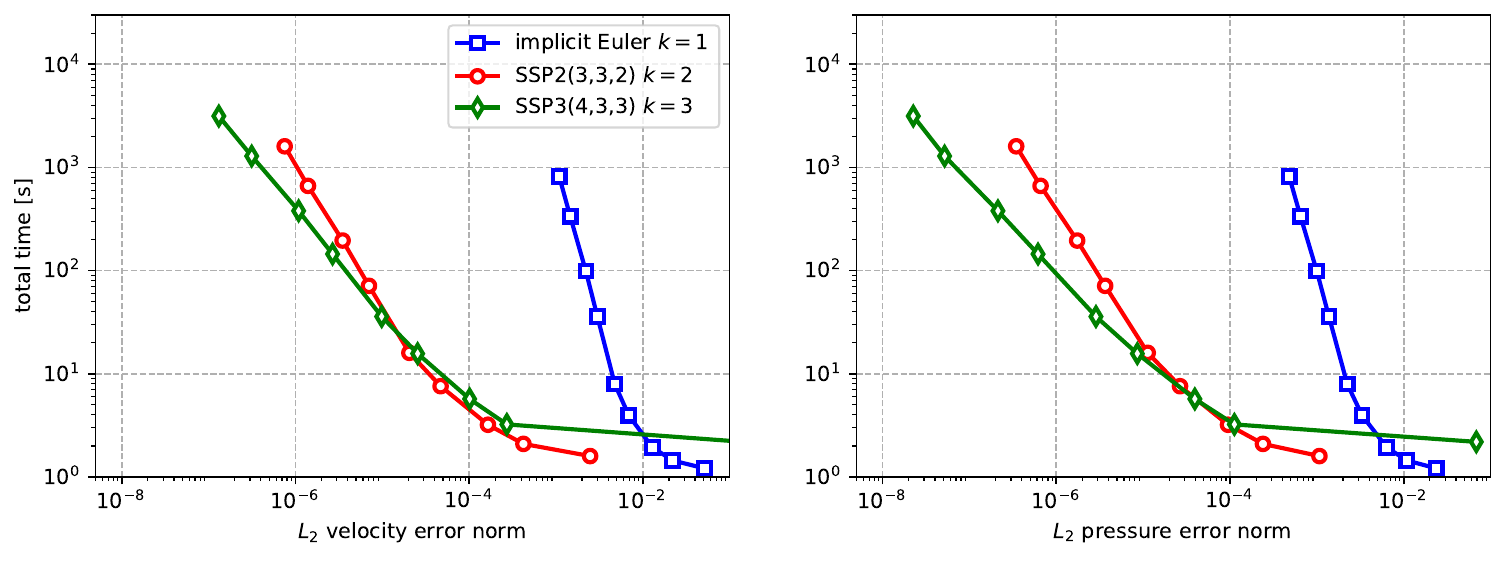}
        \caption{Total runtime for different values of the $L_2$ error in velocity (left) and pressure (right). Results are shown for the Taylor Green vortex testcase for which the analytical solution is known.}
        \label{fig:ttotal}
    \end{center}
\end{figure}

The most useful metric to quantify the overall performance of a given setup is the total runtime for a specified error. We plot this time as a function of both the $L_2$ velocity error norm and the $L_2$ pressure error norm in Fig.~\ref{fig:ttotal} for the Taylor Green vortex testcase. For a fixed error budget higher order integrators have a significantly smaller overall runtime as the bound on the error decreases.

Extrapolating the data in Fig.~\ref{fig:ttotal}, the SSP3(4,3,3) timestepper with $k=3$ is able to reduce the error to $\sim10^{-7}$ about ten times faster that the SSP2(3,3,2) method with $k=2$ while the lowest order timestepper with $k=1$ will not be able to reach this level of error in practice since it is many orders of magnitude slower. Even if only a more modest error of $\sim 10^{-4}$ is to be achieved, the higher order integrators will be about $10000\times$ faster than the lowest order method, which is clearly not competitive in this case.
\section{Conclusion and outlook}\label{sec:conclusion}
In this paper we developed a new method for the numerical solution of the time-dependent incompressible Euler equations, building on previous work in \cite{Guzman2016,Ueckermann2016}. Our approach combines a hybridisable DG discretisation with a IMEX timestepping scheme to obtain high-order convergence in space and time. Performance is achieved by using a sophisticated non-nested multigrid preconditioner for the linear system that arises in the IMEX formulation, once a projection method has been introduced to split the velocity- and pressure updates. As demonstrated by our numerical experiments, our approach is efficient in the sense that empirically the computational cost per timestep grows in proportion to the number of DG unkowns. Our method shows high order convergence, which demonstrates that it can be significantly more efficient than simpler low-order timesteppers for tight error tolerances.
\paragraph{Future work}
There are several ways in which this work can be extended. For example, no attempt has been made to optimise the IMEX scheme that was used here and it is not unlikely that other schemes (see e.g. \cite{Ascher1997,pareschi2005implicit,Weller2013}) will lead to better performance.
While our new method was applied to two well-established test cases in computational fluid dynamics, it would also be interesting to study the performance for higher-dimensional problems and more sophisticated setups. Depending on the problem, this will lead to additional complications. For example, the simulation of non-smooth flows such as the Kelvin-Helmholtz instability will require slope limiters as discussed in \cite{Ueckermann2016}. We have refrained from this here to keep the paper focused on the key aspects of the new method.

Real-live problems will likely have to be simulated on a parallel computer. It is therefore desirable to study the parallel performance of the algorithms, in particular for problems in three dimensions. Since the ILU preconditioner for the tentative velocity computation is inherently sequential, this will have to be replaced by a different approach when solving larger scale problems. The solver for the mixed pressure-velocity system, however, is likely to scale like other multigrid algorithms: parallelisation of the block-Jacobi smoother on the finest level is trivial and the scalability of the $h$-multigrid algorithm on the coarser levels has been studied extensively in the literature.

Since the focus of the paper has been the description of the new algorithm, all results reported here are empirical and based on numerical experiments. It would be highly desirable to back up the numerical results with a thorough theoretical analysis as in \cite{Guzman2016}. Clearly this is an ambitious undertaking which should be reserved for a subsequent publication.
\section*{Acknowledgements}
This work was funded by EPSRC grant EP/X019497/1 as part of the ExCALIBUR initiative. The author would like to thank the other members of the project, in particular Alexander Belozerov (Bath), as well as Colin Cotter (Imperial) and Tristan Pryer (Bath) for useful discussions and suggestions for improvement. Help and support from the Firedrake developers was crucial for the implementation of the numerical methods and is gratefully acknowledged.
\appendix

\section{Butcher tableaus}\label{sec:butcher_tableaus}
The coefficients $a_{i,j}^\impl$, $a_{i,j}^\expl$, $b_{i}^\impl$, $b_{i}^\expl$ and $c_j$ that define a particular IMEX method can be written in the form of Butcher-tableaus, see Tab.~\ref{tab:butcher_tableaus}. The Butcher-tableaus for a range of IMEX methods can be found for example in \cite{pareschi2005implicit,Weller2013}.
\begin{table}
    \begin{minipage}{0.45\linewidth}
        \begin{center}
            \begin{tabular}{c|cccccc}
                0         & 0                 &                   & \dots    &          &                     & 0               \\
                $c_1$     & $a_{1,0}^\expl$   & 0                                                                               \\
                $c_2$     & $a_{2,0}^\expl$   & $a_{2,1}^\expl$   & 0        &          &                     & \vdots          \\
                \vdots    & \vdots            &                   & $\ddots$ & $\ddots$                                         \\
                $c_{s-2}$ & $a_{s-2,0}^\expl$ & \dots             &          & $\ddots$ & 0                                     \\
                $c_{s-1}$ & $a_{s-1,0}^\expl$ & $a_{s-1,1}^\expl$ & \dots    &          & $a_{s-1,s-2}^\expl$ & 0               \\
                \hline
                          & $b_0^\expl$       & $b_1^\expl$       & \dots    &          & $b_{s-2}^\expl$     & $b_{s-1}^\expl$ \\
            \end{tabular}
        \end{center}
    \end{minipage}
    \hfill
    \begin{minipage}{0.45\linewidth}
        \begin{center}
            \begin{tabular}{|cccccc}
                0      &                   & \dots           &          &                     & 0                   \\
                0      & $a_{1,1}^\impl$                                                                            \\
                0      & $a_{2,1}^\impl$   & $a_{2,2}^\impl$ &          &                     & \vdots              \\
                \vdots &                   & $\ddots$        & $\ddots$                                             \\
                0      & \dots             &                 & $\ddots$ & $a_{s-2,s-2}^\impl$                       \\
                0      & $a_{s-1,1}^\impl$ & \dots           &          & $a_{s-1,s-2}^\impl$ & $a_{s-1,s-1}^\impl$ \\
                \hline
                0      & $b_1^\impl$       & \dots           &          & $b_{s-2}^\impl$     & $b_{s-1}^\impl$
            \end{tabular}
        \end{center}
    \end{minipage}
    \caption{Generic form of Butcher tableaus for IMEX methods considered here (recall that we assume that $a_{i,0}^\impl = b_0^\impl=0$ and hence the first column of the implicit tableau on the right-hand side contains only zeros).}
    \label{tab:butcher_tableaus}
\end{table}
For reference the Butcher tableaus for the timestepping methods that are used in our numerical experiments are shown in Tabs.~\ref{tab:butcher_tableau_imex_euler}, \ref{tab:butcher_tableau_ssp2_332} and \ref{tab:butcher_tableau_ssp3_433}.
\begin{table}
    \begin{center}
        \begin{minipage}{0.25\linewidth}
            \begin{center}
                \begin{tabular}{c|cc}
                    0 & 0 & 0 \\
                    1 & 1 & 0 \\
                    \hline
                      & 1 & 0
                \end{tabular}
            \end{center}
        \end{minipage}
        \hspace{4ex}
        \begin{minipage}{0.25\linewidth}
            \begin{center}
                \begin{tabular}{|cc}
                    0 & 0 \\
                    0 & 1 \\
                    \hline
                    0 & 1
                \end{tabular}
            \end{center}
        \end{minipage}
        \caption{Butcher tableaus for the lowest order semi-implicit IMEX Euler method}
        \label{tab:butcher_tableau_imex_euler}
    \end{center}
\end{table}
\begin{table}
    \begin{center}
        \begin{minipage}{0.25\linewidth}
            \begin{center}
                \begin{tabular}{c|ccc}
                    0              & 0              & 0              & 0              \\
                    1              & $\sfrac{1}{2}$ & 0              & 0              \\
                    $\sfrac{1}{2}$ & $\sfrac{1}{2}$ & $\sfrac{1}{2}$ & 0              \\
                    \hline
                                   & $\sfrac{1}{3}$ & $\sfrac{1}{3}$ & $\sfrac{1}{3}$
                \end{tabular}
            \end{center}
        \end{minipage}
        \hspace{4ex}
        \begin{minipage}{0.25\linewidth}
            \begin{center}
                \begin{tabular}{|ccc}
                    $\sfrac{1}{4}$ & 0              & 0              \\
                    0              & $\sfrac{1}{4}$ & 0              \\
                    $\sfrac{1}{3}$ & $\sfrac{1}{3}$ & $\sfrac{1}{3}$ \\
                    \hline
                    $\sfrac{1}{3}$ & $\sfrac{1}{3}$ & $\sfrac{1}{3}$
                \end{tabular}
            \end{center}
        \end{minipage}
        \caption{Butcher tableaus for the SSP2(3,3,2) method, see \cite{pareschi2005implicit}, \cite[Fig. 2]{Weller2013}}
        \label{tab:butcher_tableau_ssp2_332}
    \end{center}
\end{table}
\begin{table}
    \begin{center}
        \begin{minipage}{0.25\linewidth}
            \begin{center}
                \begin{tabular}{c|cccc}
                    0              & 0 & 0              & 0              & 0              \\
                    0              & 0 & 0              & 0              & 0              \\
                    1              & 0 & 1              & 0              & 0              \\
                    $\sfrac{1}{2}$ & 0 & $\sfrac{1}{4}$ & $\sfrac{1}{4}$ & 0              \\
                    \hline
                                   & 0 & $\sfrac{1}{6}$ & $\sfrac{1}{6}$ & $\sfrac{2}{3}$
                \end{tabular}
            \end{center}
        \end{minipage}
        \hspace{4ex}
        \begin{minipage}{0.25\linewidth}
            \begin{center}
                \begin{tabular}{|cccc}
                    $\alpha$  & 0              & 0              & 0              \\
                    $-\alpha$ & $\alpha$       & 0              & 0              \\
                    0         & $1-\alpha$     & $\alpha$       & 0              \\
                    $\beta$   & $\eta$         & $\delta$       & $\alpha$       \\
                    \hline
                    0         & $\sfrac{1}{6}$ & $\sfrac{1}{6}$ & $\sfrac{2}{3}$
                \end{tabular}
            \end{center}
        \end{minipage}
        \hspace{4ex}
        \begin{minipage}{0.25\linewidth}
            \begin{equation*}
                \begin{aligned}
                    \alpha & = 0.24169426078821                     \\
                    \beta  & = 0.06042356519705                     \\
                    \eta   & = 0.12915286960590                     \\
                    \delta & = \sfrac{1}{2} - \alpha - \beta - \eta
                \end{aligned}
            \end{equation*}
        \end{minipage}
        \caption{Butcher tableaus for the SSP3(4,3,3) method, see \cite{pareschi2005implicit}, \cite[Fig. 2]{Weller2013}}
        \label{tab:butcher_tableau_ssp3_433}
    \end{center}
\end{table}
\section{DG tracer advection}\label{sec:DG_tracer_advection}
Consider the advection of a passive tracer $q$ with a given velocity $Q$. In the continuum the time evolution of the tracer is given by
\begin{equation}
    \partial_t q + (Q\cdot \nabla) q = 0.\label{eqn:tracer_advection_continuum}
\end{equation}
Let $q\in V_q := \text{DG}_{h,k'}$ be a DG discretisation of the tracer field. Let $Q_{\text{CG}}=\mathcal{P}_{\text{CG}}(Q)$ be the projection of the velocity $Q\in \text{DG}_{h,k+1}$ into the continuous Galerkin space $\text{CG}_{k+1}=\text{DG}_{h,k+1}\cap C^0(\Omega)$ of the same degree. Using this projected velocity to advect the trace, the discrete form of \eqref{eqn:tracer_advection_continuum} is
\begin{equation}
    (\chi \partial_t q)_{\Omega_h} = a(\chi,q,Q_{\text{CG}}) := (q\nabla(\chi Q_{\text{CG}}))_{\Omega_h} - \sum_{F\in\mathcal{E}_h^i} \langle q^\text{up} \jump{\chi Q_{\text{CG}} \cdot n}\rangle_F
\end{equation}
for all test functions $\chi\in V_q$ where $q^\text{up}$ is the upwind value of $q$ (cf. \eqref{eqn:upwind_velocity}). Using the same IMEX timestepping method as for the velocity and pressure fields and treating the advection of the tracer explicitly, the field $q^{n+1}$ at the next timestep is given by
\begin{equation}
    (\chi q^{n+1})_{\Omega_h} = (\chi q^n) + \Delta t \sum_{i=1}^{s-1} b^\expl_i a(\chi,q_i,\mathcal{P}_{\text{CG}}(Q_i))\label{eqn:advection_update_final}
\end{equation}
The tracer value $q_i$ at the intermediate stages $i=1,\dots,s-1$ is given by
\begin{equation}
    (\chi q_i)_{\Omega_h} = (\chi q^n)_{\Omega_h} + \Delta t\sum_{j=0}^{i-1} a^\expl_{i,j} a(\chi,q_j,\mathcal{P}_{\text{CG}}(Q_j))
    \qquad\text{for all $\chi\in V_q$}\label{eqn:advection_update_stages}
\end{equation}
where $q_0 = q^n$.
\subsection{Original implicit DG method}
\begin{figure}
    \begin{center}
        \includegraphics[width=0.9\linewidth]{\figdir/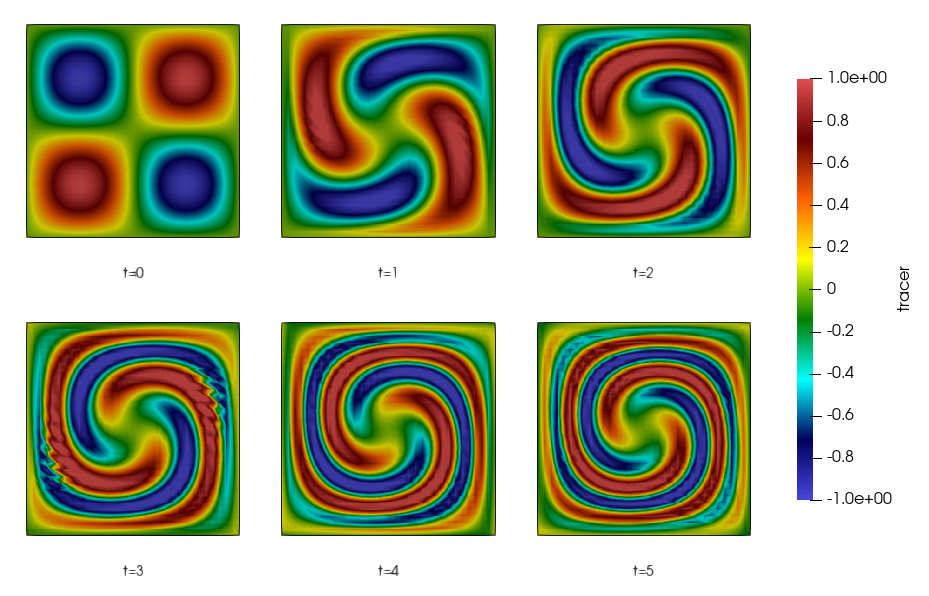}
        \caption{Advection of passive tracer in the Taylor Green vortex at different times on a $32\times 32$ grid with polynomial degree $k=1$. The results were obtained with the original, fully implicit DG method with timestep size $\Delta t=0.005$. The decay constant in the forcing function in \eqref{eqn:taylor_green_Q0_f} is set to $\kappa=0.1$.}
        \label{fig:tracer_advection_dg}
    \end{center}
\end{figure}

For the original, fully implicit DG method in \cite{Guzman2016} as described in Section \ref{sec:dg_discretisation}, the update in \eqref{eqn:advection_update_final} and \eqref{eqn:advection_update_stages} reduces to
\begin{equation}
    (\chi q^{n+1})_{\Omega_h} = (\chi q^n)_{\Omega_h} + \Delta t\;a(\chi,q^n,\mathcal{P}_{\text{CG}}(Q^n))
    \qquad\text{for all $\chi\in V_q$.}
\end{equation}
Fig.~\ref{fig:tracer_advection_dg} (which should be compared to Fig.~\ref{fig:tracer_advection}) visualises the advection of the passive tracer with this method while using the original DG timestepping scheme in \cite{Guzman2016}. A polynomial degree of $k=k'=1$ is used for both the pressure and the passive tracer, whereas velocity is of degree $k+1=2$.

\section{PETSc solver options}\label{sec:petsc_solver_options}
Tab.~\ref{tab:petsc_solver_options} lists the PETSc options that are used for the solution of the pressure system as described in Section \ref{sec:solvers}.
\begin{table}
    \begin{verbatim}
-mat_type matfree
-ksp_type preonly
-pc_type python
-pc_python_type firedrake.SCPC
-pc_sc_eliminate_fields 0, 1
-condensed_field_mat_type aij
-condensed_field_ksp_type gmres
-condensed_field_ksp_rtol 1e-12
-condensed_field_pc_type python
-condensed_field_pc_python_type firedrake.GTMGPC
-condensed_field_pc_mg_log None
-condensed_field_gt_mat_type aij
-condensed_field_gt_mg_levels_ksp_type chebyshev
-condensed_field_gt_mg_levels_pc_type python
-condensed_field_gt_mg_levels_pc_python_type firedrake.ASMStarPC
-condensed_field_gt_mg_levels_pc_star_construct_dim 1
-condensed_field_gt_mg_levels_pc_star_patch_local_type additive
-condensed_field_gt_mg_levels_pc_star_patch_sub_ksp_type preonly
-condensed_field_gt_mg_levels_pc_star_patch_sub_pc_type lu
-condensed_field_gt_mg_levels_ksp_max_it 2
-condensed_field_gt_mg_coarse_ksp_type preonly
-condensed_field_gt_mg_coarse_pc_type gamg
-condensed_field_gt_mg_coarse_pc_mg_cycles v
-condensed_field_gt_mg_coarse_mg_levels_ksp_type chebyshev
-condensed_field_gt_mg_coarse_mg_levels_ksp_max_it 2
-condensed_field_gt_mg_coarse_mg_levels_sub_pc_type sor
-condensed_field_gt_mg_coarse_mg_coarse_ksp_type chebyshev
-condensed_field_gt_mg_coarse_mg_coarse_ksp_max_it 2
-condensed_field_gt_mg_coarse_mg_coarse_sub_pc_type sor
    \end{verbatim}
    \caption{PETSc options used for the solver of the hybridised mixed pressure system.}
    \label{tab:petsc_solver_options}
\end{table}
\bibliographystyle{unsrt}

\end{document}